\documentclass[12pt]{article}
\usepackage{amsmath}
\usepackage{graphicx}
\usepackage{enumerate}
\usepackage{natbib}
\usepackage{epsfig,color,amsmath,amssymb,euscript,amsthm,epstopdf,graphicx}
\usepackage{url} 
\usepackage{xr}
\newcommand{\blind}{1}

\newtheorem{theorem}{Theorem}

\renewcommand{\hat}{\widehat}
\def\singlespace{\def\baselinestretch{1}\@normalsize}

\def\wh{\widehat}
\def\wt{\widetilde}
\def\askip{\vspace{0.1in}}

\newcommand{\cov}{{\rm Cov}}
\newcommand{\diag}{{\rm diag}}
\newcommand{\etal}{\mbox{\sl et al.\;}}
\newcommand{\ie}{\mbox{\sl i.e.\;}}

\newcommand{\var}{{\rm Var}}

\def\ga{\gamma}
\def\de{\delta}

\def\la{\lambda}

\newcommand{\bA}{{\mathbf A}}
\newcommand{\bB}{{\mathbf B}}
\newcommand{\bF}{{\mathbf F}}

\newcommand{\bG}{{\mathbf G}}
\newcommand{\bH}{{\mathbf H}}
\newcommand{\bI}{{\mathbf I}}

\newcommand{\bM}{{\mathbf M}}

\newcommand{\bP}{{\mathbf P}}
\newcommand{\bR}{{\mathbf R}}

\newcommand{\bU}{{\mathbf U}}
\newcommand{\bV}{{\mathbf V}}

\newcommand{\ba}{{\mathbf a}}
\newcommand{\bb}{{\mathbf b}}

\newcommand{\be}{{\mathbf e}}
\newcommand{\bff}{{\mathbf f}}
\newcommand{\bg}{{\mathbf g}}

\newcommand{\bx}{{\mathbf x}}
\newcommand{\by}{{\mathbf y}}
\newcommand{\bz}{{\mathbf z}}

\newcommand{\bSigma}{\boldsymbol{\Sigma}}

\newcommand{\bga}{\boldsymbol{\gamma}}
\newcommand{\bve}{\mbox{\boldmath$\varepsilon$}}

\newcommand{\bzeta} {\boldsymbol{\zeta}}

\newcommand{\bD}{{\mathbf D}}

\newcommand{\calM}{{\mathcal M}}

\def\6bullets{\bullet\bullet\bullet\bullet\bullet\bullet}

\DeclareMathAlphabet\EuScriptBF{U}{eus}{b}{n}

\def\JBES{{\sl Journal of Business \& Economic Statistics}}

\def\JASA{{\sl Journal of the American Statistical Association}}
\def\JMLR{{\sl Journal of Machine Learning Research}}

\def\AS{{\sl The Annals of Statistics}}

\def\SPA{{\sl Stochastic Processes and Their Applications}}

\def\JSPI{{\sl Journal of Statistical Planning and Inference}}

\def\JE{{\sl Journal of Econometrics}}

\newtheorem{assu}{Assumption}
\newtheorem{rmk}{Remark}

\begin{document}

\def\spacingset#1{\renewcommand{\baselinestretch}%
{#1}\small\normalsize} \spacingset{1}


\if1\blind
{
  \title{\bf Factor Modelling for Clustering High-dimensional Time Series}
  \author{Bo Zhang\thanks{
    Partially supported by \textit{National Natural Science Funds of
China No.12001517 \& 72091212, USTC Research Funds of the Double First-Class Initiative YD2040002005 and The Fundamental Research Funds for the Central Universities WK2040000026 \& WK2040000027.}}\hspace{.2cm}\\
    Department  of Statistics \&  Finance,
 International  Institute  of  Finance\\ \normalsize
School  of  Management,
University of Science and Technology of China\\ \normalsize
zhangbo890301@outlook.com
   \\
    Guangming Pan\thanks{
    Partially supported by \textit{MOE Tier 2 Grant 2018-T2-2-112 and MOE Tier 1 Grant RG133/18 at the Nanyang
Technological University, Singapore.}} \\
    School of Physical \& Mathematical Sciences,
  Nanyang Technological University\\ \normalsize
gmpan@ntu.edu.sg\\
 Qiwei Yao\thanks{
    Partially supported by \textit{EPSRC (UK) Research Grant EP/V007556/1.}}\\
 Department of Statistics, London School of Economics and Political Science\\ \normalsize
 q.yao@lse.ac.uk\\
and\\
 Wang Zhou\thanks{
    Partially supported by \textit{Grant R-155-000-211-114 at the National University of Singapore.}}\\
Department of Statistics \& Data Science,
National University of Singapore \\ \normalsize
 stazw@nus.edu.sg}
  \maketitle
} \fi

\if0\blind
{
  \bigskip
  \bigskip
  \bigskip
  \begin{center}
    {\LARGE\bf Factor Modelling for Clustering High-dimensional Time Series}
\end{center}
  \medskip
} \fi

\bigskip
\begin{abstract}
We propose a new unsupervised learning method for clustering a large number of time
series based on a latent factor structure. Each cluster is
characterized by its own cluster-specific factors in addition to some common
factors which impact on all the time series concerned. Our setting also offers
the flexibility that some time series may not belong to any clusters.
The consistency with explicit convergence rates is established for the
estimation of the common factors, the cluster-specific factors, and the latent
clusters. Numerical illustration with both simulated data as well as a real
data example is also reported. As a spin-off, the proposed new approach
also advances significantly the statistical inference for the factor model of
Lam and Yao (2012).
\end{abstract}

\noindent%
\vfill

\newpage
\spacingset{1.9} 
\section{Introduction}
\label{sec:intro}

One of the primary tasks of data mining is clustering. While most
clustering methods are originally designed for independent observations, clustering
a large number of time series gains increasing momentum (Esling and Agon 2012),
due to mining
large and complex data recorded over time in business, finance, biology,
medicine, climate, energy, environment, psychology, multimedia and other
areas (Table 1 of Aghabozorgi \etal 2015). Consequently, the literature on
time series clustering is large; see Liao (2005),
Aghabozorgi \etal (2015), Maharaj \etal (2019) and the references therein.
The basic idea is to develop some relevant
similarity or distance measures among time series first, and then
to apply the standard clustering algorithms such as hierarchical clustering
or $k$-means method. Most existing similarity/distances measures
for time series may be loosely divided into two categories: {\sl data-based}
and {\sl feature-based}. The data-based approaches define the measures
directly based on observed time series using, for example, $L_2$- or,
more general, Minkowski's distance, or various correlation measures.
Alone and Pe\~na (2019) proposed a generalized cross correlation as a similarity
measure, which takes into account cross correlation over different time lags.
 Dynamic time warping can be applied beforehand to cope with time
deformation due to, for example,
shifting holidays over different years (Keogh and Ratanamahatana, 2005).
The feature-based approaches extract relevant features from
observed time series data first, and then define similarity/distance measures based
on the extracted features. The feature extraction can be
carried out by various transformations such as Fourier,
wavelet or principal component analysis (Section 2.3 of Roelofsen, 2018).
The features from fitted time series models can also be used to
define similarity/distance measures (Yao \etal 2000, Fr\"uhwirh-Schnatter
and Kaufmann 2008).
 Attempts have also been made to define the similarity
between two time series by measuring the discrepancy between the two underlying
stochastic processes (Kakizawa \etal 1998, Khaleghi \etal 2016).
Other approaches include
Zhang (2013) which clusters time series based on the parallelism of their
trend functions, and Ando and Bai (2017) which represents the latent
clusters in terms of a factor model.
So-called `subsequence clustering' occurs frequently in the literature on time
series clustering; see Keogh and Lin (2005), and
Zolhavarieh \etal (2014). It refers to  clustering the segments from a single
long time series, which is not considered in this paper.

The goal of this study is to propose a new factor model based approach to
cluster a large number of time series
into different and unknown clusters such that the members within each cluster
share a similar dynamic structure, while the number of clusters and their sizes
are all unknown. We represent the dynamic structures
by latent common and cluster-specific factors, which are both unknown and are
identified by the difference in factor strength.
The strength of a factor is measured by the number of time series which influence
and/or are influenced by
the factor (Chamberlain and Rothschild 1983).
The common factors are strong factors 
as each of them carries the information on
most (if not all) time series concerned. The cluster-specific factors are weak
factors as they only affect the time series in a specific cluster.

Though our factor model is similar to that of Ando and Bai (2017),
our approach is radically different. First, we estimate strong factors and
all the weaker factors in the manner of one-pass, and then the latent clusters are
recovered based on the estimated weak factor loadings. Ando and Bai (2017)
adopted an iterative least squares algorithm to estimate factors/factor loadings
and latent cluster structure recursively.
Secondly, our setting allows the flexibility that some time series do not
belong to any clusters, which
is often the case in practice. Thirdly, our setting allows the dependence
between the common factors and cluster-specific factors while Ando and Bai (2017)
imposed an orthogonality condition between the two; see Remark 1(iv) in Section
\ref{sec2} below.

The methods used for estimating factors and factor loadings are adapted from
Lam and Yao (2012). Nevertheless substantial advances have been made even
within the context of Lam and Yao (2012): (i) we remove the artifact condition
that the factor loading spaces for strong and weak factors are perpendicular
with each other,  (ii) we allow weak serial correlations in idiosyncratic
components in the model, which were assumed to be vector
white noise by Lam and Yao (2012), and, more significantly, (iii) we
propose a new and {\sl consistent} ratio-based estimator for the number
of factors (see Step 1 and also  Remark 3(iii) in Section 3 below).

The rest of the paper is organized as follows. Our factor model and the relevant
conditions are presented in Section \ref{sec2}. We elaborate explicitly why
it is natural to identify the latent clusters in terms of the factor strength.
 The new clustering  algorithm is presented
in Section \ref{sec3}.
 The clustering is based on the factor loadings on all the weak factors;
 applying a $K$-means algorithm using a correlation-type similarity measure
 defined in terms of the loadings.
The asymptotic properties of the estimation for factors and factor loadings
 are collected in Section \ref{sec4}. Section \ref{sec4a} presents the consistency
the proposed factor-based clustering algorithm with error rates.
 Numerical illustration with both simulated and a real
data example is reported in Section \ref{sec5}. We also provide some comments in Section \ref{miscomment}.
All technical proofs are presented 
to 
a supplementary.

We always assume vectors in column. Let $\| \ba \|$ denote
the Euclidean norm of vector $\ba$.
For any matrix $\bG$, let $\calM(\bG)$ denote
the linear space spanned by the columns of $\bG\equiv (g_{i,j})$, $\|\bG\|$ the
square root of the largest
eigenvalue of $\bG^\top \bG$, $\|\bG\|_{\min}$  the square root of the smallest
eigenvalue of $\bG^\top \bG$,
$|\bG|$  the matrix with $|g_{i,j}|$ as its  $(i,j)$-th element.
We write $a\asymp b$ if $a=O(b)$ and $b=O(a)$.
We use $C,C_0$ to denote generic constants independent of $p$ and $n$,
which may be different at different places.
\section{Models} \label{sec2}

Let $\{\by_t\}_{1 \leq t \leq n}$ be a weakly stationary $p\times 1$ vector time series,  i.e.
$E \by_t$ is a constant independent of $t$, and all elements of
$\cov(\by_{t+k}, \by_t)$ are finite and dependent on $k$ only.  Suppose that
$\by_t$ consists of $d+1$
latent segments, \ie
\begin{equation} \label{a0}
\by_t^\top = (\by_{t,1}^\top, \cdots, \by_{t,d}^\top, \by_{t,d+1}^\top),
\end{equation}
where $\by_{t,1}, \cdots, \by_{t,d+1}$ are, respectively, $p_1, \cdots,
p_{d+1}$-vector time series with $p_1, \cdots, p_d \ge 1$, $p_{d+1} \ge 0$,
and
 $$p_1+ \cdots +p_{d} =p_0,
\qquad
p_0 + p_{d+1}=p.
$$
Furthermore, we assume  the following latent factor model with
$d$ clusters:
\begin{align} \label{a1}
\by_t & = \bA \bx_t + \Big(
{ \bB \atop \bf0}
\Big)
\bz_t + \bve_t, \\ \nonumber
\bB & = \diag(\bB_1, \cdots, \bB_d), \qquad
\bz_t^\top = (\bz_{t,1}^\top, \cdots, \bz_{t,d}^\top),
\end{align}
where $\bA$ is a $p\times r_0$ matrix with rank $r_0$, $\bx_t$ is $r_0$-vector time series representing
$r_0$ common factors and $|\var(\bx_t)| \ne 0$, $\bB_j$ is $p_j\times r_j$ matrix with rank
$r_j$, $\bz_{t,j}$ is $r_j$-vector time
series representing $r_j$ factors for $\by_{t,j}$ only and $|\var(\bz_{t,j})|\ne 0$,
$\bf0$ stands for a $p_{d+1}\times r$ matrix with all elements equal to 0,
$r= r_1 + \cdots + r_d$, and $\bve_t$ is
an idiosyncratic component in the sense of Chamberlain (1983) and Chamberlain
and Rothschild (1983) (see below). Note that in the model above, we only observe
permuted $\by_t$ (i.e. the order of components of $\by_t$ is unknown) while
all the terms on the RHS of (\ref{a1}) are unknown.

By (\ref{a1}), the $p_0$ components of $\by_t$ are grouped into $d$ clusters
$\by_{t,1}, \cdots, \by_{t,d}$, while the $p_{d+1}$ components
of $\by_{t,d+1}$ do not belong to any clusters. The $j$-th cluster $\by_{t,j}$
is characterized by the cluster-specific factor $\bz_{t,j}$,
in addition to the dependence on the common factor $\bx_t$.
The goal is to identify those $d$ latent clusters from observations
$\by_1, \cdots, \by_n$. Note that all $p_j$,
$r_j$ and $d$ are also unknown.

We always assume that the number of the common factors and the number of
cluster-specific factors for each cluster remain
bounded when the number of time series $p$ diverges. This reflects the
fact that the factor models are only appealing when the numbers of
factors are much smaller than the number of time series concerned.
Furthermore, we assume that the number of time series
in each cluster $p_i$ diverges at a lower order than $p$ and the number of clusters
$d$ diverges as well. See Assumption \ref{constantcondition} below.

\begin{assu}\label{constantcondition}
$\max_{0 \leq i \leq d}\{r_i\}<C < \infty$,
$r \asymp d=O(p^\delta)$, and $p_i \asymp p^{1-\delta}$ for $i=1, \cdots, d$,
where $C >0$ and $\delta \in (0,1)$ are constants independent of $n$ and $p$.
\end{assu}

The strength of a factor is measured by the number of time series which influence
and/or are influenced by the factor. Each component of $\bx_t$ is
a common factor. It is
related to most, if not all, components of $\by_t$ in the sense
that the most elements of the corresponding column of $\bA$ (i.e. the
factor loadings) are non-zero.
 Hence it is reasonable to assume
\begin{equation} \label{anorm}
\|\ba_j\|^2 \asymp p, \qquad j=1, \cdots,r_0,
\end{equation}
where $\ba_j$ is the $j$-th column of $\bA$.
This is in the same spirit of the definition for
the common factors by Chamberlain and Rothschild (1983).
Denoted by $\bb_i^j$ the $i$-th column of the $p_j \times r_j$ matrix $\bB_j$.
In the same vein,
we assume that
\begin{equation} \label{bnorm}
\|\bb_i^j\|^2  \asymp  p^{1-\delta}, \qquad i=1, \cdots, r_j \;\; {\rm and} \;\;
j=1, \cdots, d,
\end{equation}
as each cluster-specific factor for the $j$-th cluster is related to
most of the $p_j \asymp p^{1-\delta}$ (Assumption
\ref{constantcondition}) time series in the cluster.
Note that the factor strength can be measured by
constant $\delta \in [0, 1]$: $\delta>0$ in (\ref{bnorm}) indicates that factors
$\bz_t =(\bz_{t,1}, \cdots, \bz_{t,d})$ are weaker than factors $\bx_t$ which
corresponds to $\delta =0$; see (\ref{anorm}).

Conditions (\ref{anorm}) and (\ref{bnorm}) are imposed under the assumption
that all the factors remain unchanged as $p$ diverges,
all the entries of covariance
matrices below are bounded,
$$\mathbf{\Sigma}_x(k)=\cov(\mathbf{x}_{t+k},\mathbf{x}_t), \quad
\mathbf{\Sigma}_z(k)=\cov(\mathbf{z}_{t+k},\mathbf{z}_t),$$
$$\mathbf{\Sigma}_{x,z}(k)=\cov(\mathbf{x}_{t+k},\mathbf{z}_t), \quad
\mathbf{\Sigma}_{z,x}(k)=\cov(\mathbf{z}_{t+k},\mathbf{x}_{t}),$$
and, furthermore, $\mathbf{\Sigma}_x(k)$ and $\mathbf{\Sigma}_z(k)$ are
full-ranked for $k=0,1, \cdots, k_0$, where $k_0\ge 1$ is an integer.
Then condition (\ref{anorm}) and (\ref{bnorm}) are equivalent to
(\ref{xz1}) and (\ref{xz2}) in Assumption \ref{xandz} below after the
orthogonal normalization  to be introduced now in order to make model
(\ref{a1}) partially identifiable and operationally tractable.

In model (\ref{a1}) $\bA, \bB, \bx_t $ and $\bz_t$ are not uniquely defined,
as, for example, $(\bA, \bx_t)$ can be replaced by $(\bA \bH, \bH^{-1} \bx_t)$
for any $r_0 \times r_0$ invertible matrix $\bH$. We argue that this lack of
uniqueness gives us the flexibility to choose appropriate $\bA$ and $\bB$ to
facilitate our estimation more readily. Assumption \ref{AandB} below specifies
both $\bA$ and $\bB$ to be half-orthogonal in the sense that the columns
of $\bA$ or $\bB$
are orthonormal, which can be fulfilled by, for example, replacing
the original $(\bA, \bx_t)$ by
$(\bH, \bV\bx_t)$, where $\bA = \bH \bV$ is
a QR decomposition of $\bA$.
Even under Assumption \ref{AandB}, $\bA$ and $\bB$ are still
not unique. In fact that only the factor loading spaces $\calM(\bA), \calM(\bB_i)$
are uniquely defined by (\ref{a1}).
Hence
$\bA\bA^\top=\bA(\bA^\top \bA)^{-1}\bA^\top$, i.e. the projection
matrix onto $\calM(\bA)$, is also unique.

\begin{assu}\label{AandB}
$\mathbf{A^\top A}=\mathbf{I}_{r_0},$ $\bB_j^\top \bB_j=\mathbf{I}_{r_j}$
for $1\le j \le d$,
and it holds for a constant $q_0 \in (0, 1) $ that
\begin{equation} \label{AandBequantionq}
\|\mathbf{AA^\top}\Big( { \mathbf{B} \atop \mathbf{0}} \Big)\| \leq q_0.
\end{equation}
Furthermore for $j=1, \cdots, d$,
 ${\rm rp}(\bB_j) \{{\rm rp}(\bB_j) \}^\top$ cannot be written as a block diagonal matrix with at
least two blocks, where ${\rm rp}(\bB_j)$ denotes any row-permutation of $\bB_j$.
\end{assu}

Condition
(\ref{AandBequantionq})  implies that the columns of $(
{ \mathbf{B} \atop \mathbf{0}})$ do not fall entirely into the space
$\calM(\bA)$ as otherwise one
cannot distinguish $\bz_t$ from $\bx_t$.
It is automatically fulfilled if $\bA^\top({ \mathbf{B}
\atop \mathbf{0}}) =0$ which is a condition imposed in Lam and Yao (2012).
Finally the last  condition in Assumption \ref{AandB}
ensures that the number of clusters $d$ is uniquely defined.

\begin{assu}\label{xandz}
Let $\by_t, \bx_t$ and $ \bz_t$ be strictly stationary with the finite
fourth moments.
As $p \to \infty$, it holds for $k=0, 1, \cdots, k_0$ that
\begin{equation}\label{xz1}
 \|\mathbf{\Sigma}_x(k)\| \asymp p \asymp \|\mathbf{\Sigma}_x(k)\|_{min},
\end{equation}
\begin{equation}\label{xz2}
\|\mathbf{\Sigma}_z(k)\| \asymp p^{1-\delta}  \asymp \|\mathbf{\Sigma}_z(k)\|_{min},
\end{equation}
\begin{equation}\label{xz3}
 \|\mathbf{\Sigma}_x(k)^{-1/2}\mathbf{\Sigma}_{x,z}(k)\mathbf{\Sigma}_z(k)^{-1/2}\| \leq q_0<1, \quad \|\mathbf{\Sigma}_z(k)^{-1/2}\mathbf{\Sigma}_{z,x}(k)\mathbf{\Sigma}_x(k)^{-1/2}\| \leq q_0<1,
\end{equation}
\begin{equation}\label{xz4}
 \|\mathbf{\Sigma}_{x,z}(k)\| = O(p^{1-\delta/2}), \quad  \|\mathbf{\Sigma}_{z,x}(k)\| = O(p^{1-\delta/2}),
\end{equation}
Furthermore, $\mathbf{y}_t$ is $\psi$-mixing with the mixing coefficients satisfying $\sum_{t \geq 1}t\psi(t)^{1/2}<\infty$, and $\cov(\mathbf{x}_{t},\bve_s)=0$,
$\cov(\mathbf{z}_{t},\bve_s)=0$ for any $t$ and $s$.
\end{assu}

\begin{rmk} {\rm (i)
The factor strength is defined in terms of the orders of the factor loadings
in (\ref{anorm}) and (\ref{bnorm}). Due to the orthogonalization specified
in Assumption \ref{AandB}, they are transformed into the orders
of the covariance matrices in (\ref{xz1}) and (\ref{xz2}).
See also Remark 1 in Lam and Yao (2012).
Nevertheless, the factor strength is still measured by the constant
$\delta \in [0,1]$:
the smaller $\delta$ is, the stronger a factor is. The common factors in $\bx_t$ are
the strongest with $\delta=0$, and the cluster-specific factors in $\bz_t$ are weaker
with $\delta \in (0, 1)$.
In (\ref{a1}) $\bve_t$ represents the idiosyncratic component of $\by_t$ in the
sense that each component of $\bve_t$ only affects the corresponding component
and a few other components of $\by_t$ (i.e. $\de=1$), which
is implied by Assumptions \ref{assue} below.
Hence the strength of $\bve_t$ is the weakest.
The differences in the factor strength
make $\bx_t, \bz_t$ and $\bve_t$ on the RHS of
(\ref{a1}) (asymptotically) identifiable.
To simplify the presentation, we assume that all the components of $\bz_t$
are of the same strength (i.e. all $p_i$ are of the same order).
See the real data example in Section \ref{sec52}  for how to handle the
cluster-specific factors of different strengths.


(ii) Model (\ref{a1}) is similar to that of Ando and Bai (2017).
However, we do not require that the common factor $\bx_t$ and the
cluster-specific factor $\bz_t$ are orthogonal with each other in the
sense that ${1 \over n} \sum_{1\le t \le n} \bx_t \bz_t^\top =0$, which is imposed
by Ando and Bai (2017).
Furthermore, we allow  the idiosyncratic term $\bve_t$ to exhibit weak
autocorrelations (Assumption
\ref{assue} below),  instead of complete independence as in Ando and Bai (2017).
}
\end{rmk}

We now impose some structure assumptions on the idiosyncratic term $\bve_t$
in model (\ref{a1}).
\begin{assu}\label{assue}
Let $\bve_t = \bG \be_t$, where $\bG$ is a $p\times p$ constant matrix
with $\|\bG\|$ bounded from above by a positive constant independent of $p$.
Furthermore, one of the following two conditions holds.
\begin{description}
\item {\rm (i)} $\be_t$ is MA($\infty$), i.e.  $\be_t=\sum_{s=0}^{\infty}
\phi_s\mathbf{\eta}_{t-s}$, where $\sum_{s=0}^{\infty} |\phi_s| < \infty$,
$\mathbf{\eta}_{t}=({\eta}_{t,1},\cdots,{\eta}_{t,p})^\top$, and ${\eta}_{t,i}$
being i.i.d. across $t$ and $i$ with mean 0, variance 1 and $E({\eta}_{t,i}^4) < \infty$.
\item {\rm (ii)} $\be_t=(e_{t,1},\cdots,e_{t,p})^\top$ consists of $p$
independent weakly stationary univariate time series, $E(\be_t)=0$,
and $\min_{1 \leq i \leq p}Ee^2_{t,i}>0$. Furthermore,
$\mathbf{\tilde{e}}_{i}=(e_{1,i},\cdots,e_{n,i})$ satisfies
\begin{eqnarray}\label{subgaussiane}
\max_{\beta \geq 1, 1\leq i\leq p, \|\mathbf{a}\|=1 } \beta^{-1/2}\{\mathbb{E}|\mathbf{\tilde{e}}_i^\top\mathbf{a}|^{\beta}\}^{1/\beta} \leq C_0.
\end{eqnarray}
\end{description}
\end{assu}
\begin{rmk}
{\rm
In Assumption \ref{assue}, (i) assumes that $\be_t$ is a linear process with
the same serial correlation structure across all the components.
(ii) allows some non-linear serial dependence, the dependence structures
for different components may differ. But then the sub-Gaussian
condition (\ref{subgaussiane}) is required.
}
\end{rmk}
\section{A clustering algorithm} \label{sec3}

With available observations $\by_1, \cdots, \by_n$, we propose below
an algorithm (in five steps) to identify the latent $d$ clusters. To this end, we
introduce some notation first.
Let
$\bar \by = \frac{1}{n}\sum_{t=1}^n\by_t$,
\begin{equation}\label{Mhatdefine}
\mathbf{\hat{\Sigma}}_y(k)=\frac{1}{n}\sum_{t=1}^{n-k}(\by_{t+k} - \bar
\by)(\by_t- \bar \by)^\top, \quad
 \mathbf{\hat{M}}=\sum_{k=0}^{k_0}\mathbf{\hat{\Sigma}}_y(k)\mathbf{\hat{\Sigma}}_y(k)^\top,\
\end{equation}
where $k_0 \ge 0$ is a pre-specified integer in Assumption \ref{xandz}.

\begin{description}
\item[Step 1] (Estimation for the number of factors.)
For $0 \leq k \leq k_0$, let
$\wh \la_{k,1} \ge \cdots \ge \wh\la_{k,p}\ge 0$ be
the eigenvalues of matrix
$\mathbf{\hat{\Sigma}}_y(k)\mathbf{\hat{\Sigma}}_y(k)^\top$.
For a pre-specified positive integer $J_0 \le p$, put $\wh R_0=1$ and
\begin{equation} \label{ratios}
\wh R_j =\sum_{k=0}^{k_0}(1-{k}/{n})\wh
\la_{k,j}\Big/ \sum_{k=0}^{k_0}(1-{k}/{n})\wh \la_{k,j+1}, \quad  1 \leq j
\leq J_0.
\end{equation}
We say that $\wh R_s$ attains a local maximum if  $\wh R_s>\max\{\wh
R_{s-1},\wh R_{s+1}\}$. Let $\wh R_{\wh \tau_1}$ and $\wh R_{\wh \tau_2}$
be the two largest local maximums among $\wh R_1, \cdots, \wh R_{J_0-1}$.
The estimators for the numbers of factors are then defined as
\begin{equation} \label{ratios2}
\wh r_0 = \min\{ \wh \tau_1, \; \wh \tau_2\}, \qquad
\wh r_0 + \wh r = \max\{\wh \tau_1, \; \wh \tau_2\}.
\end{equation}

\item[Step 2] (Estimation for the loadings for common factors.)
Let $\wh\bga_1, \cdots, \wh\bga_p$ be the orthonormal eigenvectors of
matrix $\wh \bM$, arranged according to the descending order of the corresponding
eigenvalues. The estimated loading matrix for the common factors is
\begin{equation} \label{b2}
\wh \bA = (\wh \ga_1, \cdots, \wh \bga_{\wh r_0}).
\end{equation}

\item[Step 3] (Estimation for the loadings for cluster-specific factors.)
Replace $\by_t$ by $(\bI_p - \wh \bA \wh \bA^\top) \by_t$ in (\ref{Mhatdefine}), and repeat
the eigenanalysis as in Step 2 above but now denote  the corresponding
orthonormal eigenvectors by
$\wh\bzeta_1, \cdots, \wh\bzeta_p$.
The estimated loading matrix for the cluster-specific factors is
\begin{equation} \label{b3}
\wh \bB = (\wh \bzeta_1, \cdots, \wh \bzeta_{\wh r}).
\end{equation}

\item[Step 4] (Identification for the components not belonging to any clusters.)
Let $\wh \bb_1, \cdots, \wh \bb_p$ denote the row vectors of $\wh \bB$. Then
the identified index set for the components of $\by_t$ not belonging to
any clusters is
\begin{equation} \label{nocluster}
\mathfrak{\hat{J}}_{d+1}= \{ j: 1\le j \le p, \; \| \wh \bb_j \| \le \omega_p \},
\end{equation}
where $\omega_p > 0$ is a constant satisfying the conditions
$
\omega_p=o(p^{\delta/2-1/2})$,  $p^{-1/2}\omega_p^{-1}=o(1)$ and $
{p^\de n^{-1}r + p^{-\de} \over p^{1-\delta} \omega_p^2} = o(1).
$

\item[Step 5] ($K$-means clustering.)
Denote by $\hat{d}$ the number of eigenvalues of $|\wh \bB \wh
\bB^\top|$ greater than $1-\log^{-1} n$, which is taken as an upper
bound of the number of clusters.
Let $\wh p_0 = p - |\mathfrak{\hat{J}}_{d+1}|$, and
 $\wh \bF$ be the $\wh p_0 \times \wh r$ matrix obtained from $\wh \bB$
by removing the rows with their indices in $\mathfrak{\hat{J}}_{d+1}$.
Let $\wh \bff_1, \cdots, \wh \bff_{\hat p_0}$ denote the $\wh p_0$ rows of $\wh \bF$.
Let  $\wh \bR$  be the $\wh p_0 \times \wh p_0$ matrix with the
$(\ell, m)$-th element
\[
\wh \rho_{\ell, m} = \big| \wh \bff_\ell^\top \wh \bff_m\big|\big/
\big(\wh \bff_\ell^\top \wh \bff_\ell
\cdot \wh \bff_m^\top \wh \bff_m \big)^{1/2},
\qquad 1\le \ell, m \le \wh p_0.
\]
Perform the $K$-means clustering (with $L^2$-distance) for the $\wh p_0$
rows of $\wh \bR$ to form the $ d $ clusters, where $d \le \wh d$ is
chosen such that the within-cluster-sum of $L^2$-distances (to the cluster
center points) are stabilized.
\end{description}
\begin{rmk}
{\rm
(i) The ratio-based estimation in Step 1 is new. By Theorem~\ref{thma1k}
in Section \ref{sec4} below,
it holds $\wh r_0 \to r_0$ and $\wh r \to r$ in probability.
The existing approaches use the ratios of the ordered eigenvalues of matrix
$\wh \bM$ instead (Lam and Yao 2012, Chang \etal 2015, Li \etal 2017); leading
to an estimator which may not be consistent. See Example 1 below.  Note
that Lam and Yao (2012) shows that their estimator $\wt r_0$ fulfills the
relation $P(\wt r_0 \ge r_0) \to 1$ only.

(ii)  The intuition behind the estimators in (\ref{ratios}) is that the
eigenvalues $\la_{k, 1} \ge \cdots \ge \la_{k,p}( \ge 0)$ of
matrix $\bSigma_y(k) \bSigma_y(k)^\top$, where $\bSigma_y(k)= \cov(\by_{t+k}, \by_t)$,
 satisfy the conditions
\[
\la_{k,i}^{-1} =o(\la_{k,j}^{-1}) \;\;{\rm and}\;\;
\la_{k,j}^{-1} = o(\la_{k,\ell}^{-1}) \;\; {\rm for}\; 1\le i\le r_0, \;
r_0 < j \le r_0 +r \; {\rm and}\; \ell> r_0 +r.
\]
This is implied by the differences in strength among the common factor
$\bx_t$, the cluster specific factors $\bz_{t, i}$, and the
idiosyncratic components $\bve_t$; see Theorem~\ref{thma1k}.
Note that we use the ratios of the cumulative eigenvalues in (\ref{ratios}) in
order to add together the information from different lags~$k$. In practice, we set
$k_0$ to be a small integer such as $k_0\le 5$, as the significant autocorrelation
occurs typically at small lags. The results do not vary that much with respect to
the value of $k_0$ (see the simulation results in Section \ref{sec51} below).
  We truncate the
sequence $\{\wh R_j\}$ at $J_0$ to alleviate the impact of `0/0'. In
practice, we may set $J_0=p/4$ or $p/3$.

(iii) Step 3 removes the common factors first before estimating $\bB$,
as Lam and Yao (2012) showed that
weak factors can be more accurately estimated by removing strong factors from
the data first.

(iv)
 Once the numbers of strong and weak factors are correctly specified, the factor loading spaces
are relatively easier to identify. In fact $\calM(\wh \bA)$
is a consistent estimator for $\calM(\bA)$. However $\calM(\wh \bB)$ is a
consistent estimator for $\calM\{ (\bI_p - \bA \bA^\top)
({ \bB \atop{\bf0} })\}$ instead of $\calM\{({ \bB \atop{\bf0} })\}$. See
Theorem \ref{thma2} in Section \ref{sec4} below.
Furthermore the last $p_{d+1}$ rows of
$(\bI_p - \bA \bA^\top)
({ \bB \atop{\bf0} })$ are no longer 0. Nevertheless
when the elements in $\mathbf{AA}^\top$ and $\mathbf{BB}^\top$ have different orders,
those $p_{d+1}$ zero-rows can be recovered from $\wh \bB$ in Step 4.
See Theorem \ref{thma4} in Section \ref{sec4a} below.

(v)
 Given the block diagonal structure of $\bB$ in (\ref{a1}),
the $d$ clusters would be identified easily by taking the $(i,j)$-th
element of $|\bB \bB^\top|$ as the similarity measure between the $i$-th and the
$j$-th components, or by simply applying the $K$-means method to the rows of
$|\bB \bB^\top|$. However applying the $K$-means method
directly to the rows of $\bB$ will not do. Theorem \ref{thma2} and Theorem \ref{thma4old} indicate that the block
diagonal structure, though masked by asymptotically diminishing `noise',
 still presents in $\wh \bB \wh \bB^\top$ via a latent row-permutation of
$\wh \bB$. Accordingly
the cluster analysis in Step 5 is based on the absolute values of
the correlation-type measures
among the rows of $\wh \bF \wh \bF^\top$ which is an estimator for $\bB \bB^\top$.

(vi) In Step 5, we search for the number of clusters $d$ by the `elbow method'
which is the most frequently used method in $K$-means clustering. Nevertheless
$\wh d$ provides an upper bound for $d$; see Theorem \ref{consistdthm} below.
Our empirical experiences indicate that $\wh d =d$ holds often especially
when $r_j$, $1\le j \le d$, are small. See Tables \ref{t4} and \ref{t4n800} in
Section \ref{sec51} below. Note that $\bB \bB^\top$ is a block diagonal matrix
with $d$ blocks and all the non-zero eigenvalues equal to 1. Therefore
the dominant eigenvalue for each of the latent $d$ blocks in
$\wh \bB \wh \bB^\top$ is greater than or at least very close to 1.
Moreover, by Perron-Frobenius's theorem, the largest eigenvalue of
$|\bB_j\bB_j^\top|$, i.e. the so-called Perron-Frobenius eigenvalue,
 is strictly greater than the other eigenvalues of $|\bB_j\bB_j^\top|$
under the last condition in Assumption \ref{AandB}.
This is the intuition behind the definition of $\wh d$.

}
\end{rmk}

\noindent
{\bf Example 1}. Consider a simple model of the form (\ref{a1}) in which
$\bve_t \equiv 0$, $r_0=1, \, r=2$, $\bA^\top\Big({
\mathbf{B} \atop \mathbf{0}}\Big)=\bf0$ and
\begin{align*}
& x_t = p^{1/2}(u_{1,t} + a_1 u_{1, t-1}
+ u_{2,t} + a_2 u_{2, t-1}), \\
& z_{1,t} =p^{1/2-\de/2}(u_{2,t}+ a_2 u_{2,t-1}), \quad
z_{2,t}= p^{1/2-\de/2}(u_{3,t} + a_3 u_{3,t-1}),
\end{align*}
where $a_1 , a_2, a_3$ are constants, and $u_{i,t}$, for different $i,t$, are
independent and $N(0,1)$.
Let $\bM = \sum_{0\le k \le 1} \bSigma_y(k) \bSigma_y(k)^\top$, and $\la_1 \ge \la_2 \ge \la_3$
be the three largest eigenvalues of $\bM$.
It can be shown that  $\la_1 \asymp p^2$, $\la_3 = p^{2-2\de}\{ (1+a_3^2)^2 + a_3^2\}$
and $\la_2 \asymp p^{2-\delta}$ provided
$(a_1 - a_2)^2(1- a_1 a_2) \ne 0$.  Hence $\la_1/\la_2 \asymp \la_2 / \la_3 \asymp p^\de$.
This shows that $r_0 (=1)$ or $r (=2)$ cannot be estimated stably based on the ratios of the
eigenvalues of $\wh \bM$ for this example. In fact, let  two $p \times 3$ matrices $\bU_0$ and $\bU_1$ be the eigenvectors of $\bSigma_y(0)\bSigma_y(0)^\top$ and $\bSigma_y(1)\bSigma_y(1)^\top$. When $(a_1 - a_2)^2(1- a_1 a_2) \ne 0$,  $\bU_0$ and $\bU_1$ are different while $\bU_0\bU_0^\top=\bU_1\bU_1^\top$.


\section{Asymptotic properties on estimation for factors} \label{sec4}

Theorem \ref{thma1} and Remark \ref{remake4} below show that in the absence of weak factor
$\bz_t$, the estimation for the strong factor loading
space $\calM(\bA)$ achieves root-$n$ convergence rate in spite of diverging $p$.
Since only the factor loading space $\calM(\bA)$
is uniquely defined by (\ref{a1}) (see the discussion
below Assumption \ref{AandB}), we measure the estimation error in terms of its (unique)
projection matrix $\bA \bA^\top$.

\begin{theorem}\label{thma1}
Let Assumptions \ref{constantcondition}-\ref{assue} hold.
Let $p, n\to \infty$,  $n=O(p)$ and $p^{\delta}=o(n)$.
Then it holds that
\begin{equation}\label{Ahatjieguo2}
 \big\|\mathbf{\hat{A}}\mathbf{\hat{A}}^\top-\mathbf{A}\mathbf{A}^\top\big\|=O_p(n^{-1/2}+p^{-\delta/2}).
\end{equation}
\end{theorem}

Assumption \ref{AandB} ensures that  the rank of matrix $\bB_{*} \equiv
(\bI_p - \bA \bA^\top)\Big({
\mathbf{B} \atop \mathbf{0}}\Big)$ is $r$. Denote by $\bP_{{\tiny A_{\bot}B}} = \bB_{*}
(\bB_{*}^\top \bB_{*})^{-1} \bB_{*}^\top$ the projection matrix
onto $\calM\big\{(\bI_p - \bA \bA^\top)\Big({
\mathbf{B} \atop \mathbf{0}}\Big)\big\}$ of which $\calM(\wh \bB) $ is a consistent estimator,
see Theorem \ref{thma2} below, and also Remark 3(iv).

\begin{theorem}\label{thma2}
Let Assumptions \ref{constantcondition}-\ref{assue} hold.
Let $p, n\to \infty$,  $n=O(p)$, $p^{\delta}r=o(n)$.
Then it holds that
\begin{equation}\label{Bhatjieguo2}
 \big\|\mathbf{\hat{B}}\mathbf{\hat{B}}^\top- \bP_{{\tiny A_{\bot}B}}\big\|
=O_p(p^{\delta/2}n^{-1/2}+p^{-\delta/2})
\end{equation}
and
\begin{equation}\label{Bhatjieguo2a}
 \big\|\mathbf{\hat{B}}\mathbf{\hat{B}}^\top- \bP_{{\tiny A_{\bot}B}}\big\|_F
=O_p(p^{\delta/2}n^{-1/2}r^{1/2}+p^{-\delta/2}).
\end{equation}
\end{theorem}

Theorem \ref{thma1k} below specifies the asymptotic behavior for
the ratios
of the cumulated eigenvalues used in estimating
the numbers of factors in Step 1 in Section \ref{sec3} above.
%
It implies
that $\wh r_0 \to r_0, \; \wh r \to r$ in
probability provided that $J_0 > r_0 + r$ is fixed.

\begin{theorem}\label{thma1k}
Let Assumptions \ref{constantcondition}-\ref{assue} hold.
Let $p, n\to \infty$,  $n=O(p)$, $p^{\delta}r=o(n)$.
For $\wh R_j$ defined in (\ref{ratios}), it holds for some constant $C>0$ that
\begin{align}\label{Ahatjieguo1kratio1a}
& \lim_{n,p \rightarrow \infty}P(\wh R_j < C) =1 \quad {\rm for}\;\; j=1, \cdots, r_0-1,\\
\label{Ahatjieguo2kratio2a}
& \;\;\; \wh R_{r_0}^{\, -1}
=O_p(p^{-2\delta}),  \quad
\wh R_{r_0+r}^{\, -1}
=O_p\big(n^{-2}p^{2\delta}\big),\\
\label{Ahatjieguo1kratio1b}
& \lim_{n,p \rightarrow \infty}P( \wh R_j
< C)=1 \quad {\rm for} \;\; j=r_0+1, \cdots, r_0+r-1, \;\; {\rm and}\\
\label{Ahatjieguo1kratio1c}
& \;\;\;
\wh R_j
=O_p(1) \quad {\rm for}\;\; j= r_0+r+1, \cdots, r_0+r+s,
\end{align}
where $s$ is a positive fixed integer.
\end{theorem}
\begin{rmk}\label{remake4}
It is worth pointing out that
 the block diagonal structure of $\bB$ is not required for Theorems
\ref{thma1}--\ref{thma1k}.
On the other hand, if $\bA^\top\Big({
\mathbf{B} \atop \mathbf{0}}\Big) =0$ and $\{\bx_t\}$ and $\{\bz_t\}$ are independent,
the term $p^{-\delta/2}$ on the RHS of (\ref{Ahatjieguo2})-(\ref{Bhatjieguo2a})
disappears.
\end{rmk}

%
%
%

\section{Asymptotic properties on clustering}\label{sec4a}
\begin{assu}\label{assuaifix}
The elements of $\mathbf{AA^\top }$ are of the order $O(p^{-1})$,
and $\|\mathbf{b}_i\|^2 \asymp p^{\delta-1}$ for $1 \leq i \leq p_0$,
where $\mathbf{b}_i$ denotes the $i$-th row of matrix $\bB$.
\end{assu}
The orthogonalization $\bA^\top \bA= \bI_{r_0}$ implies that the
average of the squared elements of $\bA$ is $O(p^{-1})$. Since $r_0$ is finite,
it is reasonable to assume that the elements of $\mathbf{AA^\top }$ are $O(p^{-1})$.
As $\bB$ is a block-diagonal matrix with blocks $\{ \bB_j \}$ and $\bB_j^\top
\bB_j=\bI_{r_j}$, the squared elements of $\bB_j$ are of the
order $p_j^{-1}\asymp p^{-(1-\delta)}$ in average.   As $r_j$ is bounded, it is
reasonable to assume $\|\mathbf{b}_i\|^2 \asymp p^{\delta-1}$.
Assumption \ref{assuaifix} ensures that
$(\bI_p - \bA \bA^\top)\Big( { \mathbf{B} \atop \mathbf{0}}\Big)$ is
asymptotically a block diagonal matrix; see Theorem \ref{thma4old} below.
This enables to recover the block diagonal structure of $\bB$ based on
$\wh \bB$ which provides a consistent estimator for the space
$\calM\big\{(\bI_p - \bA \bA^\top)\Big( { \mathbf{B} \atop \mathbf{0}}\Big)\big\}$
(Theorem \ref{thma2} above), and also to separate the components of $\by_t$
not belonging to any clusters. See Theorems \ref{thma4} -- \ref{consistkmean} below.
\begin{theorem}\label{thma4old}
Let Assumptions \ref{constantcondition} -- \ref{AandB} and \ref{assuaifix} hold.
Divide matrix $\Big({ \mathbf{B} \atop \mathbf{0}}\Big)(\mathbf{B}^\top,
\mathbf{0})- \bP_{{\tiny A_{\bot}B}}$ into $(d+1) \times (d+1)$ blocks, and
denote its $(i,j)$-th block of the size $p_i \times p_j$ by $\mathbf{V}_{i,j}$.
Then as $p\to \infty$,
$\|\mathbf{V}_{i,j}\|_F=O_p(p^{-1}p_i^{1/2}p_j^{1/2})= O_p(p^{-\delta})$.
\end{theorem}

\begin{theorem}\label{thma4}
Let the conditions of Theorem \ref{thma2}
 and Assumption \ref{assuaifix}  hold.
%
For
$\mathfrak{\hat{J}}_{d+1}$ defined in (\ref{nocluster}),
it holds that
\begin{equation}\label{assuaifix6}
\frac{|\mathfrak{J}_{d+1}^c\cap\mathfrak{\hat{J}}_{d+1}|}{p^{1-\delta}}=O_p\Big(n^{-1}rp^{\delta}+p^{-\delta}\Big), 
\quad {\rm and}
\end{equation}
\begin{equation}\label{assuaifix7}
\frac{|\mathfrak{J}_{d+1}\cap\mathfrak{\hat{J}}_{d+1}|}{|\mathfrak{J}_{d+1}|}=1+O_p\Big(\frac{p^{\delta}rn^{-1}+p^{-\delta}}{p_{d+1}\omega_p^2}\Big),
\end{equation}
where $\omega_p$ is given in (\ref{nocluster}).
Furthermore,
$\frac{p^{\delta}rn^{-1}+p^{-\delta}}{p_{d+1}\omega_p^2}=o(1)$ provided
that $p^{1-\delta}/p_{d+1}=O(1)$.

\end{theorem}

\begin{theorem}\label{consistdthm}
Let the conditions of Theorem \ref{thma4} hold, and $p^{\delta}r\log^{2} n=o(n)$.
Then $ P(\hat{d}\geq d) \to 1, $ as $n,p \to \infty$.

\end{theorem}
\begin{rmk}
Theorem \ref{thma4} shows that most the components belonging to the
$d$ clusters will not be classified as not belonging to any
clusters (see (\ref{assuaifix6})).  Furthermore most the components
not belonging to any clusters will be correctly identified (see
(\ref{assuaifix7})).
 Theorem \ref{consistdthm} shows that the probability of under-estimating $d$ converges to 0.
\end{rmk}


To investigate the errors in the $K$-means clustering, let
$\mathbf{R} = ( |r_{\ell, m}|)$ be the $p_0 \times p_0$ matrix with
 \[
 r_{\ell, m} =   \bb_\ell^\top  \bb_m\big/ \big( \bb_\ell^\top  \bb_\ell
\cdot  \bb_m^\top  \bb_m \big)^{1/2}, \qquad 1\le \ell, m \le p_0.
\]
We assume that $d$ is known.
Let $\mathfrak{O}_d$ be the set consisting of all $p_0 \times p_0$
matrices with $d$ distinct rows. Put
\begin{eqnarray}\label{kmeancontroldeff0}
\mathbf{D}_0=\arg\min_{\mathbf{D} \in \mathfrak{O}_{d}}\|\mathbf{R}-\mathbf{D}\|_F^2.
\end{eqnarray}
For any $p_0$-vector $\bg$ with its elements taking integer values between
1 and $d$, let
\begin{eqnarray*}
\mathfrak{O}_d(\bg) = \{ \bD \in \mathfrak{O}_d: &
{\rm two\; rows\; of\; {\bf D}\; are\; the\; same\; if\; and \;only\; if\;
the \; corresponding} \\ &
\mbox{two elements of {\bf g} are the same} \}.
\end{eqnarray*}

Note that the $d$ distinct rows of $\bD_0$ would be the centers of the $d$ clusters
identified by the $K$-means method based on the rows of $\bR$. However
$\bR$ is unknown, we identify the $d$ clusters based on its estimator $\wh \bR$;
see Step 5 of the algorithm in Section \ref{sec3}.
The clustering based on $\wh \bR$
can only be successful if that based on $\bR$ is successful, i.e. $\bD_0 \in
\mathfrak{O}_d(\bg_0)$, where $\bg_0$ is the $p_0$-vector with 1 as its first
$p_1$ elements, 2 as its next $p_2$ elements, $\cdots$, and $d$ as its last $p_d$
elements.
Given the block diagonal structure of $\bB$, condition $\bD_0 \in
\mathfrak{O}_d(\bg_0)$
is likely to hold.

For any $p_0$-vector $\bg$ with its elements taking integer values between
1 and $d$, it partitions $\{ 1, \cdots, p_0\}$ into $d$ subset. Let $\tau(\bg)$
denote the number of misclassified components by partition $\bg$.

\begin{assu}\label{kmeancorrassu}
$\bD_0 \in
\mathfrak{O}_d(\bg_0)$. For some constant $c>0$,
\[
\min_{\mathbf{D} \in \mathfrak{O}_{d
}(\mathbf{g})}\|\mathbf{R}-\mathbf{D}\|_F^2 \geq
\|\mathbf{R}-\mathbf{D}_0\|_F^2+c\tau(\bg) p^{1-\delta}.
\]
\end{assu}

\begin{theorem}\label{consistkmean}
Let the conditions of Theorem \ref{thma4} and Assumption \ref{kmeancorrassu} hold.
The number of clusters $d$ is assumed to be known. Denoted by $\wh
\tau$ the number of misclassified components of $\by_t$ by the $K$-means
clustering in Step 5.  Then as $n, p \to \infty$,
\begin{eqnarray}\label{misclusterrate0}
\wh \tau/p = O_p\big( p^{-\delta/2} \big).
\end{eqnarray}
\end{theorem}
\begin{rmk}\label{remake7a}
Theorem \ref{consistkmean} implies that the misclassification
rate of the $K$-means method converges to 0, though the convergence
rate is slow (see (\ref{misclusterrate0})). However a faster rate is attained
when $\bA^\top\Big({
\mathbf{B} \atop \mathbf{0}}\Big) =0$, and $\{\bx_t\}$ and $\{\bz_t\}$ are independent, as then $\wh \tau/p = O_p\big(n^{-1/2}r^{1/2} \big)$. See also Remark \ref{remake4}.
Assumption \ref{kmeancorrassu}
requires that $\|\mathbf{R}-\mathbf{D}\|_F^2$ increases as the number of misplaced members of this  partition $\tau(\bg)$ becomes large,
which is necessary for the K-means method. 
\end{rmk}

\section{Numerical properties} \label{sec5}

\subsection{Simulation } \label{sec51}
We illustrate the proposed methodology through a simulation study with
model (\ref{a1}).
We draw the elements of $\bA$ and $\bB_j$ independently from $U(-1, 1)$.  All component series of $\bx_t$ and
$\bz_t$ are independent and AR(1) and MA(1), respectively, with Gaussian innovations.
All components of $\bve_t$ are independent MA(1) with $N(0, 0.25)$ innovations.
All the AR and the MA coefficients are drawn randomly from $U\{(-0.95,
-0.4)\cup (0.4, 0.95)\}$.
The standard deviations of the components of $\bx_t$ and
$\bz_t$ are drawn randomly from $U(1, 2)$.

We consider following two scenarios with $r_0=r_1 = \cdots = r_d=2$ and
$p_1=\cdots = p_d$:
\begin{itemize}
\item Scenario I: $n=400$, $d=5$ and $p_{d+1}=p_1$. Hence $r=10$ and $p=6p_1$.
\item Scenario II: $n=800$, $d=10$ and $p_{d+1}=5 p_1$. Hence $r=20$ and $p=15p_1$.
\end{itemize}
The numbers of factors $r_0$ and $r$ are estimated based on the
ratios $\wh R_j$ in (\ref{ratios2}) with $k_0=1, \cdots, 5$ and $J_0=[p/4]$ in
(\ref{ratios}).
For the comparison purpose, we also report the estimates based on the
ratios of eigenvalues of $\wh \bM$ in (\ref{Mhatdefine}) also with $k_0=0, \cdots, 5$,
which is the standard
method used in literature and is defined as in (\ref{ratios2}) but now
with $\wh R_j= \wt \la_j /\wt \la_{j+1}$ instead, where
$\wt \la_1 \ge \cdots \ge \wt \la_p \ge 0$ are the eigenvalues of $\wh \bM$.
 See, e.g. Lam and Yao (2012).
For each setting, we replicate the experiment 1000 times.

The relative
frequencies of $\wh r_0 = r_0$ and $\wh r_0 + \wh r=r_0+ r$ are reported
in Tables~\ref{t1}-\ref{t1800}.
Overall the method based on the ratios of the cumulative eigenvalues $\wh R_j$
provides accurate and robust performance and is not sensitive to the
choice of $k_0$.
The estimation based on the eigenvalues of $\wh \bM$ with $k\ge 1$
is competitive for $r_0$, but is considerably poorer for $r_0+r$
in Scenario II.
Using $\wh \bM$ with $k=0$ leads to weaker estimates for $r_0$ in Scenario I.

It is noticeable that the performance of the estimation for the number
of common factor $r_0$ in Scenario II is better than Scenario I. This is due to the fact the difference in the factor strength between
the common factor $\bx_t$ and the cluster-based factor $\bz_t$ in Scenario II is larger than Scenarios I. In contrast, the results hardly change with different values
of $p_1$.

\begin{table}
\caption{\scriptsize{The relative frequencies of $\wh r_0=r_0$ and $\wh r_0 + \wh r=r_0+r$ in a simulation for  Scenario I
with 1000 replications, where $\wh r_0$ and $\wh r$ are estimated by
the $\wh R_j$-based method (\ref{ratios2}),
and the ratios of the eigenvalues of $\hat \bM$. } } \label{t1}
\begin{center}
\scriptsize{\begin{tabular}{l|llll|llll}
Estimation& \multicolumn{4}{c|}{$\wh r_0=r_0$} & \multicolumn{4}{c}{$\wh r_0
+ \wh  r=r_0+ r$}\\
method    & $p_1=25$ & $p_1=50$ & $p_1=75$ &  $p_1=100$
&$p_1=25$ & $p_1=50$ & $p_1=75$ &  $p_1=100$   \\
\hline
$\wh R_j$ $(k_0=1)$
& .742 &.785  &.787 &.808
& 1  &      1 &   1 &1\\
$\wh R_j$ $(k_0=2)
$& .762 &.792 &.799 &.823
& .999&    1 &   1 &1\\
$\wh R_j$ $(k_0=3)$
& .766 &.787 &.793 &.823
& .999&    1 &   1 &1\\
$\wh R_j$ $(k_0=4)$
&  .753 &.783 &.790 &.821
& .998 &   1 &   1 &1\\
$\wh R_j$ $(k_0=5)$
& .751 &.779 &.782 &.816
& .998 &   1 &   1 &1\\
{ $\wh \bM$ $(k_0=0)$}
& .707 &.749  &.761& .758
& 1    &    1  &  1 &1\\
{$\wh \bM$ $(k_0=1)$}
& .754 &.791 &.790 &.788
& .985 &.993 &.989 &.971\\
{ $\wh \bM$ $(k_0=2)$}
& .774 &.802 &.807   &.826
&  .977 &.993 &.998 &.998\\
{ $\wh \bM$ $(k_0=3)$}
& .783 &.801 &.807 &.829
& .953 &.989 &.992 &.996
\\
{ $\wh \bM$ $(k_0=4)$}
& .786 &.799 &.805 &.830
& .929 &.979 &.982 &.991\\
{ $\wh \bM$ $(k_0=5)$}
& .781 &.797 &.804 &.828
& .916 &.965 &.976 &.979
\\
\end{tabular}}
\end{center}
\end{table}


\begin{table}
\caption{\scriptsize{The relative frequencies of $\wh r_0=r_0$ and $\wh r_0 + \wh r=r_0+r$ in a simulation for  Scenario II
with 1000 replications, where $\wh r_0$ and $\wh r$ are estimated by
the $\wh R_j$-based method (\ref{ratios2}),
and the ratios of the eigenvalues
of $\wh \bM$.}} \label{t1800}
\begin{center}
\scriptsize{\begin{tabular}{l|llll|llll}
Estimation& \multicolumn{4}{c|}{$\wh r_0=r_0$} & \multicolumn{4}{c}{$\wh r_0
+ \wh  r=r_0+ r$}\\
method   & $p_1=25$ & $p_1=50$ & $p_1=75$ &  $p_1=100$
&$p_1=25$ & $p_1=50$ & $p_1=75$ &  $p_1=100$   \\
\hline
$\wh R_j$ $(k_0=1)$
& .985 &.989 &.996 &.993
 & 1  &      1 &   1 &1\\
$\wh R_j$ $(k_0=2)$
& .985 &.984 &.992 &.990
& 1&    1 &   1 &1\\
$\wh R_j$ $(k_0=3)$
& .983 &.981 &.989 &.988
& .999&    1 &   1 &1\\
$\wh R_j$ $(k_0=4)$
&  .980 &.974 &.986 &.985
& .998 &   1 &   1 &1\\
$\wh R_j$ $(k_0=5)$
& .976 &.973 &.984 &.985
& .998 &   1 &   1 &1\\
{ $\wh \bM$ $(k_0=0)$}
& .976 &.991 &.998& .993
& 1    &    1  &  1 &1\\
{$\wh \bM$ $(k_0=1)$}
& .985 &.990 &.996 &.993
& .638  &.647 &.642 &.646\\
{ $\wh \bM$ $(k_0=2)$}
& .986 &.985 &.992   &.990
&  .705 &.780 &.801 &.818\\
{ $\wh \bM$ $(k_0=3)$}
& .984 &.981 &.989 &.989
& .638 &.734 &.762 &.771
\\
{ $\wh \bM$ $(k_0=4)$}
& .980 &.976 &.985 &.988
& .597 &.689 &.722 &.723\\
{ $\wh \bM$ $(k_0=5)$}
& .978 &.973 &.985  &.987
& .569 &.660 &.694 &.702
\\
\end{tabular}}
\end{center}
\end{table}

\begin{table}
\caption{\scriptsize{The means and standard deviations (in parentheses) of
$\|\mathbf{\hat{A}\hat{A}^\top}-\mathbf{AA^\top}\|_F$ and
$\|\mathbf{\hat{B}\hat{B}^\top}-\bP_{A\bot B}\|_F$ in a simulation for  Scenario I
with 1000 replications, where $\wh \bA$ is estimated by the
eigenvectors of $\wh \bM$ in (\ref{Mhatdefine}) (with $k_0=1, \cdots, 5$),
 or by those of $\wh \bSigma_y(k) \wh \bSigma_y(k)^\top$ (for $k=0, 1, \cdots, 5$),
and $\wh \bB$ is estimated in the similar manner.
Both $r_0$ and $r$ are assumed to be known.}} \label{t2}

\askip

\resizebox{\textwidth}{!}{
\scriptsize{\begin{tabular}{l|llll|llll}
Estimation&
\multicolumn{4}{c|}{\small $\|{\hat{\bA}\hat{\bA}^\top}-\mathbf{AA^\top}\|_F$}
& \multicolumn{4}{c}{\small $\|{\hat{\bB}\hat{\bB}^\top}-\bP_{A\bot B}\|_F $}\\
method    & $p_1=25$ & $p_1=50$ & $p_1=75$ &  $p_1=100$
&$p_1=25$ & $p_1=50$ & $p_1=75$ &  $p_1=100$   \\
\hline
{ $\wh \bM$ $(k_0=1)$}
&.237(.101) &.219(.077)  &.212(.077) &.205(.069)
&.508(.061)& .491(.043)& .486(.043)& .482(.039)\\
{ $\wh \bM$ $(k_0=2)$}
& .229(.092) &.215(.072) &.210(.073) &.203(.065)
&.511(.056)& .495(.042)&.491(.042)&.487(.038)\\
{ $\wh \bM$ $(k_0=3)$}
& .228(.089) &.215(.071) &.211(.072) &.204(.064)
&.517(.055)& .501(.042)&.498(.042)& .493(.038)\\
{ $\wh \bM$ $(k_0=4)$}
& .229(.089) &.216(.071) &.212(.072) &.205(.064)
&.523(.055)& .507(.043)& .503(.042)& .499(.039)\\
{ $\wh \bM$ $(k_0=5)$}
 & .230(.088) &.217(.070) &.214(.073) &.206(.064)
 &.528(.055)& .511(.043)& .508(.043)& .503(.039)\\
{ $\hat{\bSigma}_y(0)\hat{\bSigma}_y(0)^\top$}
&.264(.122) &.234(.088) &.222(.086) &.213(.077)
& .501(.077)& .478(.050)& .470(.049)& .466(.044)\\
{ $ \hat{\bSigma}_y(1)\hat{\bSigma}_y(1)^\top$}
&.292(.117) &.283(.101) &.281(.105) &.274(.096)
& .994(.083)& .976(.075)&.973(.072)& .969(.072)\\
{ $ \hat{\bSigma}_y(2)\hat{\bSigma}_y(2)^\top$}
& .414(.187) &.411(.181) &.414(.181) &.405(.180)
& 2.27(.139)& 2.26(.139)& 2.26(.131)& 2.25(.139)\\
{ $ \hat{\bSigma}_y(3)\hat{\bSigma}_y(3)^\top$}
& .621(.312) &.614(.303) &.632(.309) &.610(.306)
& 2.31(.154)& 2.30(.152)& 2.30(.152)& 2.29(.150)\\
{ $ \hat{\bSigma}_y(4)\hat{\bSigma}_y(4)^\top$}
& .775(.352) &.767(.336) &.800(.361)&.764(.348)
& 2.35(.163)& 2.34(.158)& 2.34(.161)& 2.33(.158)\\
{ $ \hat{\bSigma}_y(5)\hat{\bSigma}_y(5)^\top$}
&.885(.359) &.881(.347) &.894(.358) &.868(.350)
& 2.39(.169)& 2.37(.158)& 2.36(.164)& 2.36(.164)\\
\end{tabular}}
}
\end{table}

\begin{table}
\caption{ \scriptsize{The means and standard deviations (in parentheses) of
$\|\mathbf{\hat{A}\hat{A}^\top}-\mathbf{AA^\top}\|_F$ and
$\|\mathbf{\hat{B}\hat{B}^\top}-\bP_{A\bot B}\|_F$ in a simulation for  Scenario II
with 1000 replications, where $\wh \bA$ is estimated by the
eigenvectors of $\wh \bM$ in (\ref{Mhatdefine}) (with $k_0=1, \cdots, 5$),
 or by those of $\wh \bSigma_y(k) \wh \bSigma_y(k)^\top$ (for $k=0, 1, \cdots, 5$),
and $\wh \bB$ is estimated in the similar manner.
Both $r_0$ and $r$ are assumed to be known.}} \label{t2n800}

\askip

\resizebox{\textwidth}{!}{
\scriptsize{\begin{tabular}{l|llll|llll}
Estimation&
\multicolumn{4}{c|}{\small $\|{\hat{\bA}\hat{\bA}^\top}-\mathbf{AA^\top}\|_F$}
& \multicolumn{4}{c}{\small $\|{\hat{\bB}\hat{\bB}^\top}-\bP_{A\bot B}\|_F $}\\
method    & $p_1=25$ & $p_1=50$ & $p_1=75$ &  $p_1=100$
&$p_1=25$ & $p_1=50$ & $p_1=75$ &  $p_1=100$   \\
\hline
{ $\wh \bM$ $(k_0=1)$}
&.118(.027) &.116(.027)  &.115(.026) &.113(.025)
&.720(.040)&.705(.034)& .701(.033)& .699(.032)\\
{ $\wh \bM$ $(k_0=2)$}
& .119(.027) &.117(.027) &.117(.026) &.115(.025)
&.735(.043)& .718(.035)& .714(.035)&.711(.033)\\
{ $\wh \bM$ $(k_0=3)$}
& .121(.027) &.120(.027) &.119(.026) &.117(.025)
&.748(.045)& .729(.037)& .725(.036)& .722(.034)\\
{ $\wh \bM$ $(k_0=4)$}
& .123(.028) &.122(.027) &.121(.027) &.119(.025)
&.758(.047)& .739(.038)& .734(.037)& .731(.035)\\
{ $\wh \bM$ $(k_0=5)$}
& .124(.028) &.123(.027) &.122(.027) &.120(.025)
&.767(.048)& .746(.038)& .741(.038)& .738(.036)\\
{ $\hat{\bSigma}_y(0)\hat{\bSigma}_y(0)^\top$}
&.119(.029) &.114(.026) &.113(.025) &.111(.024)
& .677(.036)& .665(.031)& .662(.031)& .661(.030)\\
{ $ \hat{\bSigma}_y(1)\hat{\bSigma}_y(1)^\top$}
&.169(.048) &.168(.049) &.166(.046) &.165(.047)
& 1.48(.080)& 1.45(.073)& 1.44(.071)& 1.44(.069)\\
{ $ \hat{\bSigma}_y(2)\hat{\bSigma}_y(2)^\top$}
& .277(.121) &.277(.122) &.274(.114) &.276(.125)
& 3.29(.101)& 3.28(.099)& 3.27(.098)& 3.27(.103)\\
{ $ \hat{\bSigma}_y(3)\hat{\bSigma}_y(3)^\top$}
& .456(.249) &.455(.249) &.450(.240) &.444(.247)
& 3.32(.106)& 3.31(.105)& 3.30(.106)& 3.30(.105)\\
{ $ \hat{\bSigma}_y(4)\hat{\bSigma}_y(4)^\top$}
& .626(.329) &.621(.314) &.610(.306)&.604(.309)
& 3.35(.119)& 3.33(.112)& 3.32(.114)& 3.32(.112)\\
{ $ \hat{\bSigma}_y(5)\hat{\bSigma}_y(5)^\top$}
&.729(.337) &.730(.330) &.725(.337) &.723(.338)
& 3.37(.123)& 3.36(.114)& 3.34(.118)& 3.35(.115)\\
\end{tabular}}
}
\end{table}
Recall $\bP_{A\bot B}$ is the projection matrix onto the space
$\calM\big\{ (\bI_p - \bA \bA^\top) \big( {\bB \atop{\bf0}} \big) \big\}$; see
Theorem \ref{thma2} and also Remark 3(iv).
Tables \ref{t2}-\ref{t2n800} contain the means and standard deviations of the estimation
errors for the factor loading spaces $\|{\hat{\bA}\hat{\bA}^\top}-\mathbf{AA^\top}\|_F$
and $\|{\hat{\bB}\hat{\bB}^\top}-\bP_{A\bot B}\|_F $, where $\wh \bA$ is
estimated by the eigenvectors of the     matrix $\wh \bM$ in (\ref{Mhatdefine})
with $k_0=1, \cdots, 5$, see Step 2 of the algorithm stated in Section~\ref{sec3}.
See also Step 3 there for the similar procedure in estimating $\bB$.
For the comparison purpose, we also include the estimates obtained with $\wh \bM$
replaced by $\hat{\bSigma}_y(k)\hat{\bSigma}_y(k)^\top$ with $k=0, 1, \cdots, 5$.
Tables \ref{t2}-\ref{t2n800} show clearly that the estimation based on $\wh \bM$ is accurate
and robust with respect to the different values of $k_0$. Furthermore
 using a single-lagged
covariance matrix for estimating factor loading spaces is not recommendable.
The error $\|{\hat{\bA}\hat{\bA}^\top}-\mathbf{AA^\top}\|_F$ in Scenario
I is larger than the error in Scenarios II. This is due to  the larger
sample size $n$ in Scenario II. See Theorem \ref{thma1}. In contrast,
Theorem \ref{thma2} shows that the error rate of
$\|{\hat{\bB}\hat{\bB}^\top}-\bP_{A\bot B}\|$ contains the term
$p^{\delta/2}n^{-1/2}$. While $n$ is larger  in Scenario II, so is $p^{\delta}$.
This explains why the error
$\|{\hat{\bB}\hat{\bB}^\top}-\bP_{A\bot B}\|_F$ in Scenario II is also larger
than that in Scenario I.

In the sequel, we only report the results with $\wh r_0$ and $\wh r$ estimated
by (\ref{ratios2}), and the factor loading spaces estimated by
the eigenvectors of $\wh\bM$. We always set $k_0=5$.
We examine now the effectiveness of Step 4 of the algorithm.
Note that the indices of the components not belonging to any clusters are
identified as those in $\wh{\mathfrak J}_{d+1}$ in (\ref{nocluster}),
which is defined in terms of a threshold $\omega_p=o(p^{\delta/2-1/2})$.
We experiment with the three choices of this tuning parameter,
namely $\omega_{p1}=(\hat{r}/p)^{1/2}/\ln p $,
$\omega_{p2}=\{\hat{r}/(p\ln p)\}^{1/2}$ and
$\omega_{p3}=\{\hat{r}/(p\ln \ln p)\}^{1/2}$.
Recall ${\mathfrak J}_{d+1}^c$ contains all the indices of the components of $\by_t$ belonging
to one of the $d$ clusters.
The means and standard deviations of the two types of misclassification errors
$ E_1=|{\mathfrak J}_{d+1}^c \cap \wh {\mathfrak J}_{d+1}|/|{\mathfrak J}^c_{d+1}|$
and
$E_2=|{\mathfrak J}_{d+1} \cap \wh {\mathfrak J}_{d+1}^c|/|{\mathfrak J}_{d+1}|$
over the 1000 replications are reported in Tables~\ref{t3}-\ref{t3n800}.
Among the three choices, $\omega_{p2}$ appears
to work best as the two types of errors are both small. The increase in the errors
due to the estimation for $r_0$ and $r$ is not significant.

\begin{table}
\caption{\scriptsize{The means and standard deviations (in parentheses) of the error
rates $E_1=|{\mathfrak J}_{d+1}^c \cap \wh {\mathfrak J}_{d+1}|/|{\mathfrak J}^c_{d+1}|$
and $E_2=|{\mathfrak J}_{d+1} \cap \wh {\mathfrak J}_{d+1}^c|/|{\mathfrak J}_{d+1}|$
in a simulation for  Scenario I with 1000 replications
with the 3 possible choices of threshold $\omega_p$ in (\ref{nocluster}),
and the numbers of factors $r_0$ and $r$ either known or to be estimated.
}} \label{t3}

\askip

\resizebox{\textwidth}{!}{
\begin{tabular}{ll|llll|llll}
&&
\multicolumn{4}{c|}{$r_0$ and $r$ are known}
& \multicolumn{4}{c}{$r_0$ and $r$ are estimated}\\
& & $p_1=25$ & $p_1=50$ & $p_1=75$ &  $p_1=100$
&$p_1=25$ & $p_1=50$ & $p_1=75$ &  $p_1=100$   \\
\hline
$\omega_{p1}$ & {\small $E_1$ }
&.007(.007) &.006(.005) &.005(.003) &.005(.003)                     & .006(.007) &.005(.006) &.005(.004) &.004(.007)\\
$\omega_{p2}$ & & .073(.021)&.067(.014) &.064(.012) &.061(.010)   & .067(.024) &.062(.022) &.058(.016) &.057(.023)\\
$\omega_{p3}$ & & .259(.031) &.244(.021) &.237(.017) &.231(.015)      & .256(.033) &.240(.026) &.232(.018) &.227(.023)\\
\hline
$\omega_{p1}$ & {\small $E_2$}
& .279(.152) &.216(.119) &.201(.107) &.194(.097)            & .380  (.247) &.328(.255) &.321(.264) &.297(.254) \\
$\omega_{p2}$ & & .000(.003) &2e-5(.0006) &1.3e-5(.0004) &.000(.002)    &.050(.096) &.051(.101) &.057(.111) &.050(.108)  \\
$\omega_{p3}$ & & .000(.000) &.000(.000) &.000(.000) &.000(.000)         & .000(.001) &.000(.000) &.000(.000) &.000(.000) \\
\end{tabular}
}
\end{table}
\begin{table}
\caption{\scriptsize{The means and standard deviations (in parentheses) of the error
rates $E_1=|{\mathfrak J}_{d+1}^c \cap \wh {\mathfrak J}_{d+1}|/|{\mathfrak J}^c_{d+1}|$
and $E_2=|{\mathfrak J}_{d+1} \cap \wh {\mathfrak J}_{d+1}^c|/|{\mathfrak J}_{d+1}|$
in a simulation for  Scenario II with 1000 replications
with the 3 possible choices of threshold $\omega_p$ in (\ref{nocluster}),
and the numbers of factors $r_0$ and $r$ either known or to be estimated.
}} \label{t3n800}

\askip

\resizebox{\textwidth}{!}{
\begin{tabular}{ll|llll|llll}
&&
\multicolumn{4}{c|}{$r_0$ and $r$ are known}
& \multicolumn{4}{c}{$r_0$ and $r$ are estimated}\\
& & $p_1=25$ & $p_1=50$ & $p_1=75$ &  $p_1=100$
&$p_1=25$ & $p_1=50$ & $p_1=75$ &  $p_1=100$   \\
\hline
$\omega_{p1}$ & {\small $E_1$ }
                &.003(.004) &.003(.002) &.002(.002) &.002(.001)       & .003(.004) &.002(.002) &.002(.002) &.002(.001)\\
$\omega_{p2}$ & & .049(.013)&.046(.009) &.044(.007) &.042(.006)       & .049(.013) &.046(.009) &.044(.007) &.042(.007)\\
$\omega_{p3}$ & & .185(.019) &.179(.014) &.174(.012) &.171(.010)      & .185(.019) &.179(.015) &.174(.012) &.171(.010)\\
\hline
$\omega_{p1}$ & {\small $E_2$}
                & .208(.083) &.202(.069) &.226(.072) &.252(.070)      &.221(.118) &.218(.119) &.235(.104) &.261(.100) \\
$\omega_{p2}$ & & .000(.003) &.000(.000) &.000(.000) &.000(.000)      &.001(.005)  &.001(.009) &.001(.009) &.001(.011)  \\
$\omega_{p3}$ & & .000(.000) &.000(.000) &.000(.000) &.000(.000)      & .000(.000) &.000(.000) &.000(.000) &.000(.000) \\
\end{tabular}
}
\end{table}


In the sequel, we only report the results with
$\omega_{p2}=\{\hat{r}/(p\ln p)\}^{1/2} $. In Step 5, we estimate
$\wh d$ as an upper bound for $d$.
As $r_j=2$ for $j=1, \cdots, d$, $\wh d =d$ occurs almost surely in our simulation.
See Tables \ref{t4}-\ref{t4n800}.
 Then the $\wh d$ clusters are obtained by performing the $K$-means clustering
for the $\wh p_0$ rows of $\wh \bR$, where $\wh p_0= p- |\mathfrak{\hat
J}_{d+1}|$. As the error rates in estimating $\mathfrak{J}_{d+1}^c$
has already been reported in Tables \ref{t3}-\ref{t3n800}, we concentrate
on the components of $\by_t$ with indices in $\mathfrak{\hat{J}}_{d+1}^c
\cap \mathfrak{J}_{d+1}^c$ now, and count the number of them which were misplaced
by the $K$-means clustering, i.e. $\wh \tau$.
Both the means and the standard deviations of the error rates
$\wh \tau/|\mathfrak{\hat{J}}_{d+1}^c \cap \mathfrak{J}_{d+1}^c|$
over 1000 replications are reported in Tables \ref{t4}-\ref{t4n800}. We also report the relative frequencies of $\wh d=d$.
Tables \ref{t4}-\ref{t4n800} show
clearly that the $K$-means clustering identifies the latent clusters
very accurately, and the difference in performance due to the estimating $(r_0, r)$
is also small.

\begin{table}
\caption{\scriptsize{The means and standard deviations (STD)  of the error
rates $\wh \tau/|\mathfrak{\hat{J}}_{d+1}^c \cap \mathfrak{J}_{d+1}^c|$ and the relative frequencies of $\wh d=d$
in a simulation for  Scenario I with 1000 replications
with the numbers of factors $r_0$ and $r$ either known or to be estimated.
}} \label{t4}

\askip

\makebox[\textwidth]{
\begin{tabular}{l|llll|llll}
&
\multicolumn{4}{c|}{$r_0$ and $r$ are known}
& \multicolumn{4}{c}{$r_0$ and $r$ are estimated}\\
  & $p_1=25$ & $p_1=50$ & $p_1=75$ &  $p_1=100$
&$p_1=25$ & $p_1=50$ & $p_1=75$ &  $p_1=100$   \\
\hline
mean
          &1.7e-5 &0 &2.8e-6 &0         & .0037 &.0029&.0037 &.0034\\
STD  & .0003&0 &8.9e-5&0               & .0134 &.0072 &.0101&.0167\\
$\hat{d}=d$ &1&1&1&1                  &1&.999&1&.999
\end{tabular}
}
\end{table}

\begin{table}
\caption{\scriptsize{The means and standard deviations (STD)  of the error
rates $\wh \tau/|\mathfrak{\hat{J}}_{d+1}^c \cap \mathfrak{J}_{d+1}^c|$ and the relative frequencies of $\wh d=d$
in a simulation for  Scenario II with 1000 replications
with the numbers of factors $r_0$ and $r$ either known or to be estimated.
} }\label{t4n800}

\askip

\makebox[\textwidth]{
\begin{tabular}{l|llll|llll}
&
\multicolumn{4}{c|}{$r_0$ and $r$ are known}
& \multicolumn{4}{c}{$r_0$ and $r$ are estimated}\\
  & $p_1=25$ & $p_1=50$ & $p_1=75$ &  $p_1=100$
&$p_1=25$ & $p_1=50$ & $p_1=75$ &  $p_1=100$   \\
\hline
mean
          &0 &0 &0 &0         & 8e-6 &.0002&.0003 &.0005\\
STD  & .0003&0 &0&0                & .0001 &.0036 &.0046&.0075\\
$\hat{d}=d$ &1&1&1&1                  &1&1&1&1
\end{tabular}
}
\end{table}

\subsection{Real data illustration} \label{sec52}

We consider the daily returns of the stocks listed in S\&P500 in
31 December 2014 -- 31 December 2019. By removing those which were not traded on
every trading day during the period, there are $p=477$ stocks which were traded
on $n=1259$ trading days. Those stocks are from 11 industry sectors:
\begin{center}
\begin{tabular}{lll}
1. Communication Services \quad& 2. Consumer Discretionary \quad& 3. Consumer Staples\\
4. Energy& 5. Financials&  6. Health Care\\
7. Industrials& 8. Information Technology& 9. Materials\\
10.Real Estate& 11.Utilities&
\end{tabular}
\end{center}
The conventional wisdom suggests that the companies in the same industry sector
share some common features. We apply the proposed 5-step algorithm in
Section \ref{sec3} to the return series
to cluster those 477 stocks into different groups.

\begin{figure}
\centering
\epsfig{figure=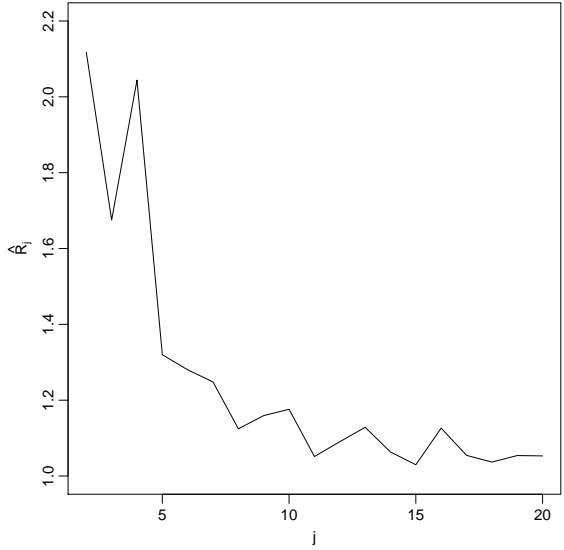,height=8cm,width=8cm}
\caption{\footnotesize{Plot of $\hat{R}_j$ against $j$ for $2 \leq j \leq 20$.} }
\label{FRDwhole1}
\end{figure}

Step 1 is to
estimate the numbers of strong factors and cluster-specific weak factors.
To this end, we
calculate $\wh R_j$ as in (\ref{ratios}) with $k_0=5$.
It turns out
$\hat{R}_1=32.53$ is much larger than all the others, while $\hat{R}_j$ for $j \geq 2$
are plotted in Figure \ref{FRDwhole1}. By (\ref{ratios2}), $\wh r_0=1$ and $\wh r_0+
\wh r=4$. Note that the estimates for $\wh r_0$ and $\wh r_0+
\wh r$ are unchanged with $k_0=1, \cdots, 4$.
While  the existence of
$\wh r_0=1$ strong and common factor is reasonable, it is most unlikely that
there are merely $\wh r =3$ cluster-specific weak factors.
Note that estimators in
(\ref{ratios2}) are derived under the
assumption that all the $r$ cluster-specific (i.e. weak) factors are of the same
factor strength; see Remark 1(ii) in Section 2 above. In practice
weak factors may have different degrees of strength; implying that we should also
take into account the 3rd, the 4th, the 5th largest local maximum of $\wh R_j$.
Hence we take $\wh r_0+
\wh r=16$ (or perhaps also 10 or 13), as Figure \ref{FRDwhole1} suggests that there are 3
factors with factor strength $\delta_1 >0$, and further 12 factors with
strength $\delta_2 \in (\delta_1, 1)$.

\begin{figure}
\centering
\epsfig{figure=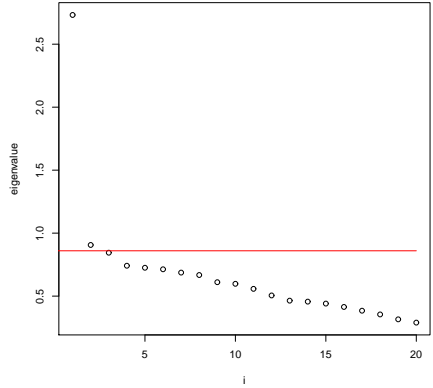,height=6cm,width=6cm}
\caption{\footnotesize{The eigenvalues of $|\mathbf{\hat{B}\hat{B}}^\top|$ when $\hat{r}_0=1$ and $\hat{r}=15$. The red line is $1-\log^{-1}n$.}}
\label{FRDwhole2}
\end{figure}

With $\wh r_0=1$ and $\wh r=15$, we proceed to Steps 2 \& 3 of Section
\ref{sec3} and obtain the estimator
$\wh \bB$ as in (\ref{b3}).  Setting $\omega_{p}=\big\{\wh r/(p \ln p)\big\}^{1/2}$,
$|\mathfrak{\hat{J}}_{d+1}|=11$, i.e. 11 stocks do not appear to belong to any clusters,
where $\mathfrak{\hat{J}}_{d+1}$ is defined as in (\ref{nocluster}) in Step 4.
Leaving those 11 stocks out, we perform Step 5, i.e the $K$-means
clustering for the $\wh p_0= 477-11 = 466$ rows of matrix $\wh \bR$.
From Figure \ref{FRDwhole2} we choose $\hat{d}=3$.
But we also consider $\hat{d}=9$ and $\hat{d}=11$ as two more examples.

To present the identified $d$ clusters, we define $11\times d$ matrix
with $n_{ij}/n_i$ as its $(i,j)$-th element, where $n_i$ is the number
of the stocks in the $i$-th industry sector, and $n_{ij}$ is the number
of the stocks in the $i$-th industry sector which are allocated in the
$j$-th cluster. Thus $n_{ij}/n_i \in [0,1]$ and $\sum_j n_{ij}/n_i=1$.

\begin{figure}
\epsfig{figure=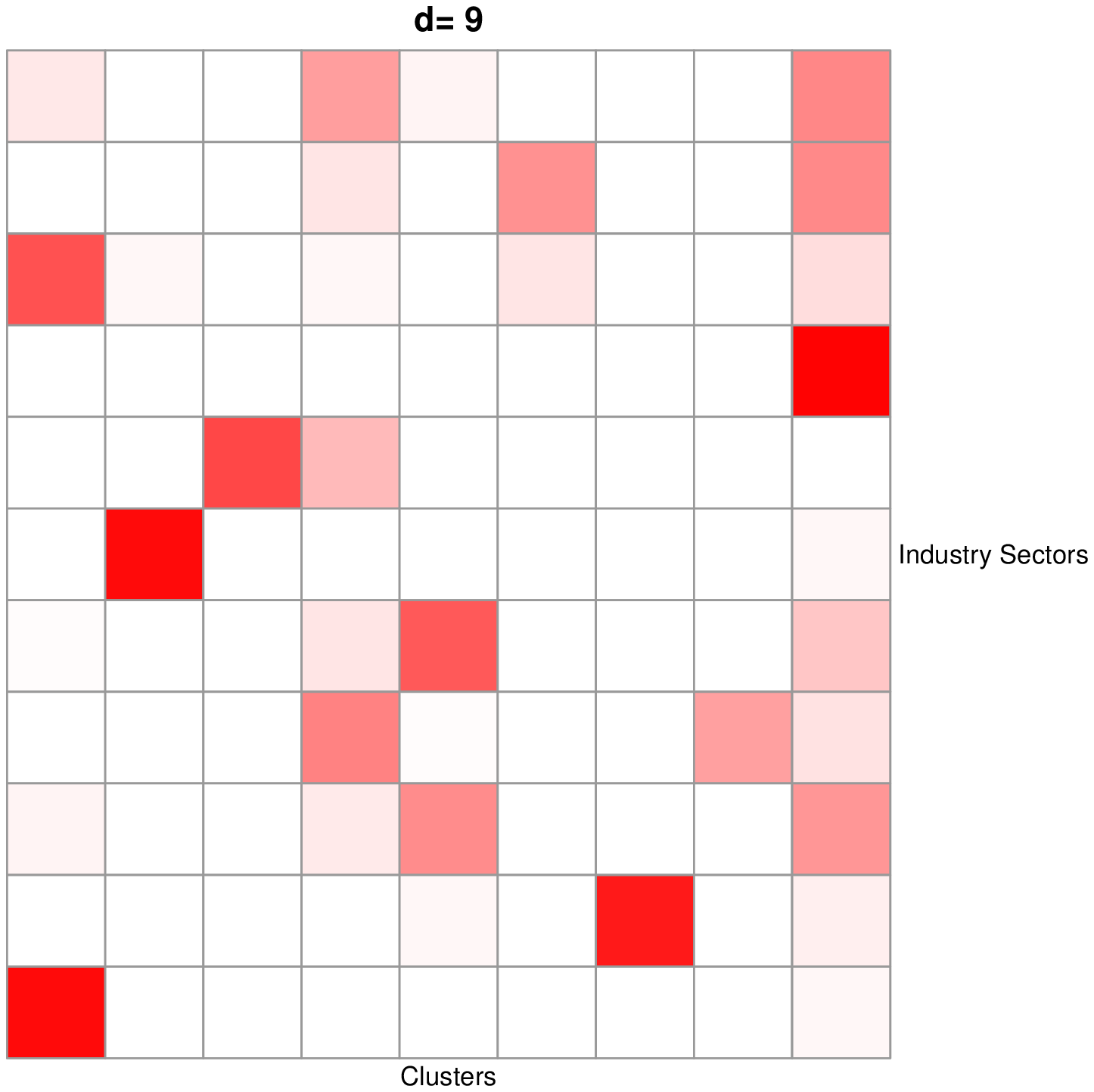, width=5cm, height=6cm} \;\;\;\;
\epsfig{figure=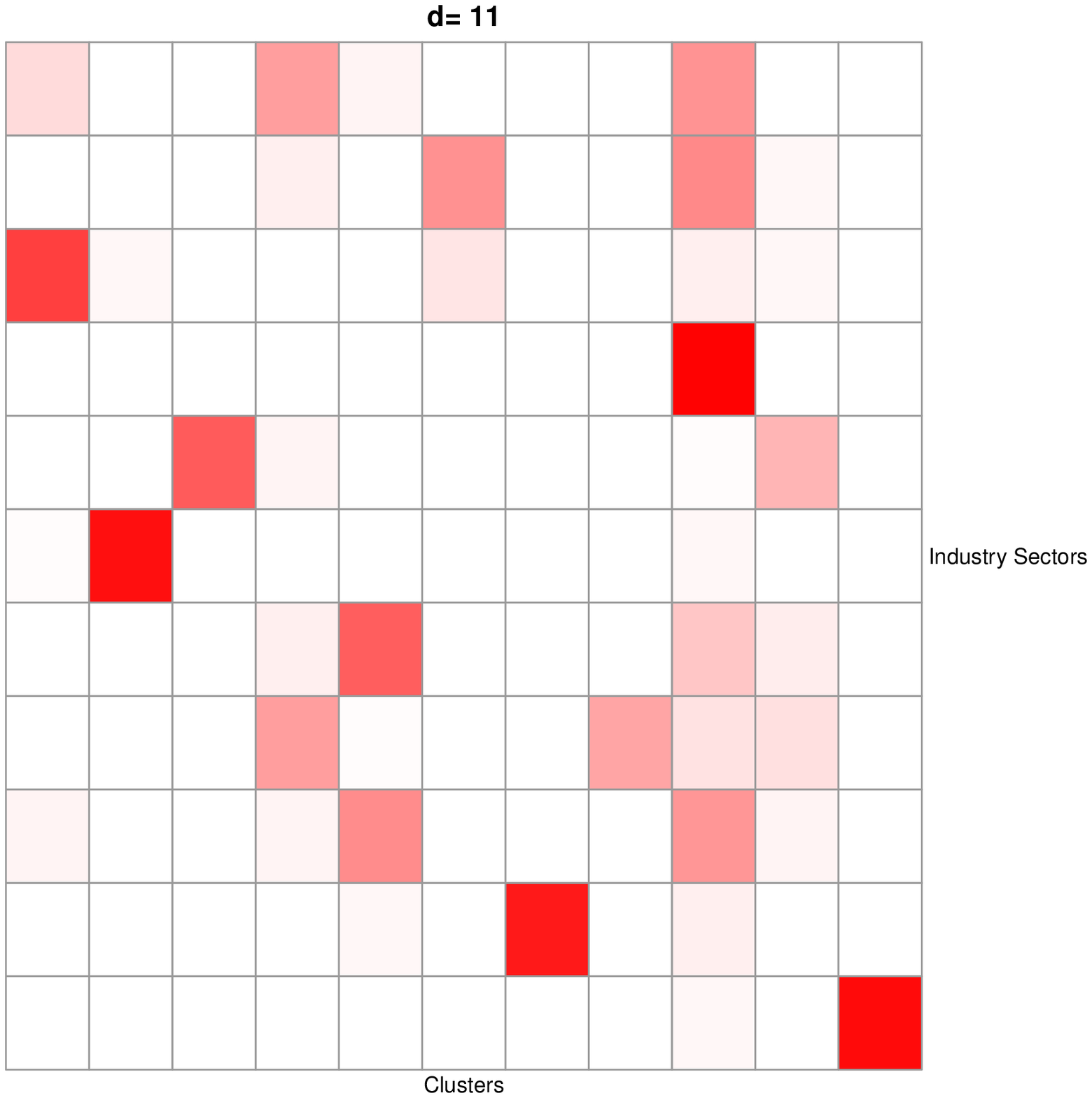, width=6cm, height=6cm}\;\;\;\;
\epsfig{figure=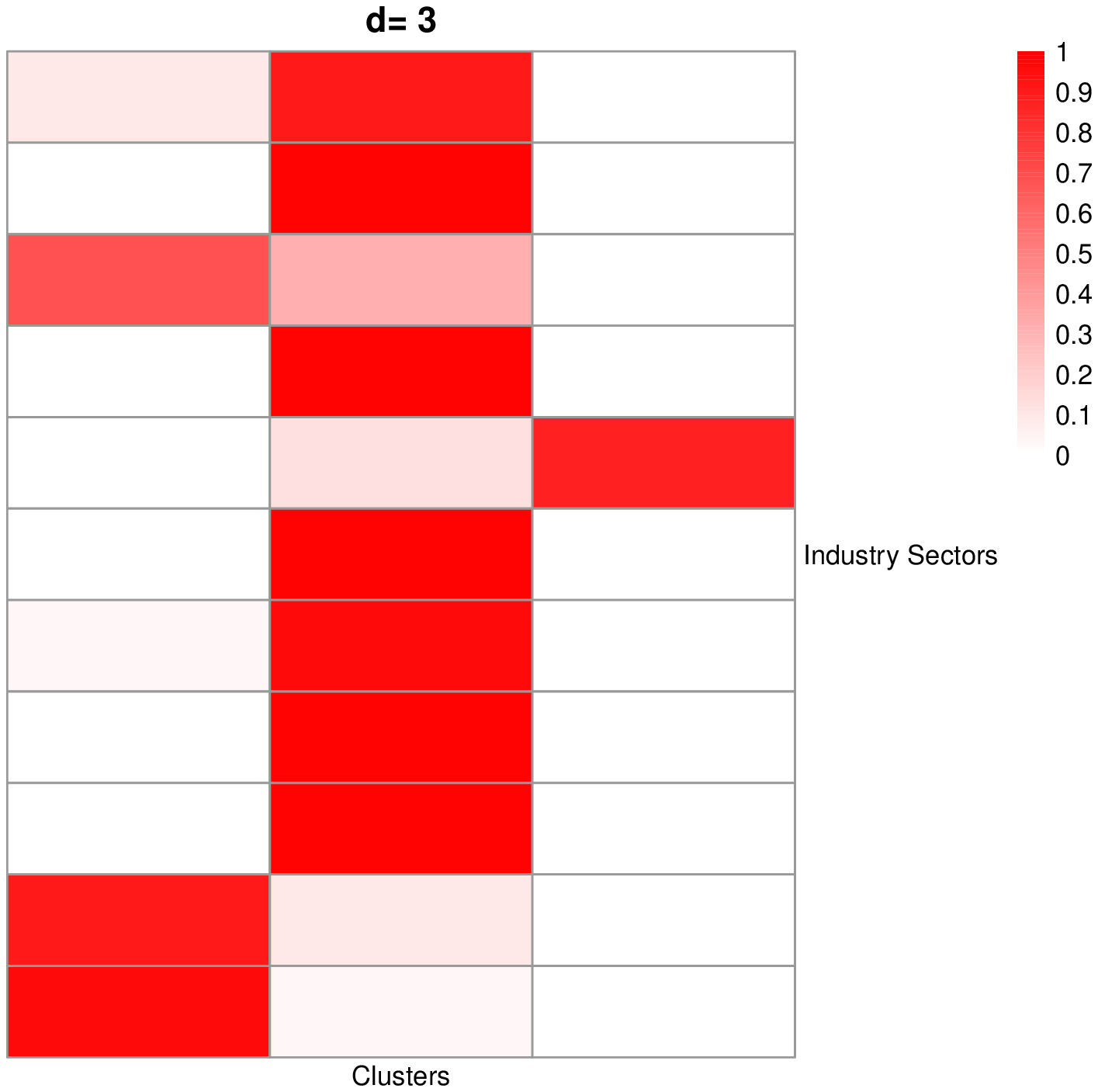, width=4cm, height=6cm}
\caption{\footnotesize{Heat-maps of the distributions of the stocks in each of
the 11 industry sectors (corresponding to 11 rows) over $d$ clusters (corresponding
to $d$ columns), with $d=9, 11 $ and 3. The estimated
numbers of the common and cluster-specific factors are,
respectively, $\wh r_0=1$ and $\wh r =15$.}} \label{fig5r15}
\end{figure}

The heat-maps of this $11\times d$ matrix for $d=\wh d=9$ is presented in Figure \ref{fig5r15}.
The first cluster mainly consists of the companies in Consumer Staples,
 and Utilities, Clusters 2 and 3 contain the companies in, respectively,
Health Care and Financials, Cluster 4 contains mainly some companies in Communication
Service and Information Technology, Cluster 5 consists of the companies in Industrials
and Materials, Cluster 6 are mainly the companies in Consumer Discretionary,
Cluster 7 are mainly the companies in Real Estate.   Cluster 8 is mainly the companies from Information
Technology, Cluster 9 contains almost all companies in Energy and a small number of companies from  each of 5 or
6 different sectors.
To examine how stable the clustering is, we also include the results for
$d=11$ and $d=3$ in Figure \ref{fig5r15}.
When $d$ is increased from 9 to 11, the original Cluster 1 is divided into
new Clusters 1 and 11 with the former consisting of Consumer Staples, and the latter being Utilities. Furthermore the
original Cluster 4 splits into new Clusters 4 and 10, while the other
7 original clusters are hardly changed.
With $d=3$, most companies in each of the 11 sectors stay in one cluster. For example, most companies in Financials are always in a separate group.

\begin{figure}
\epsfig{figure=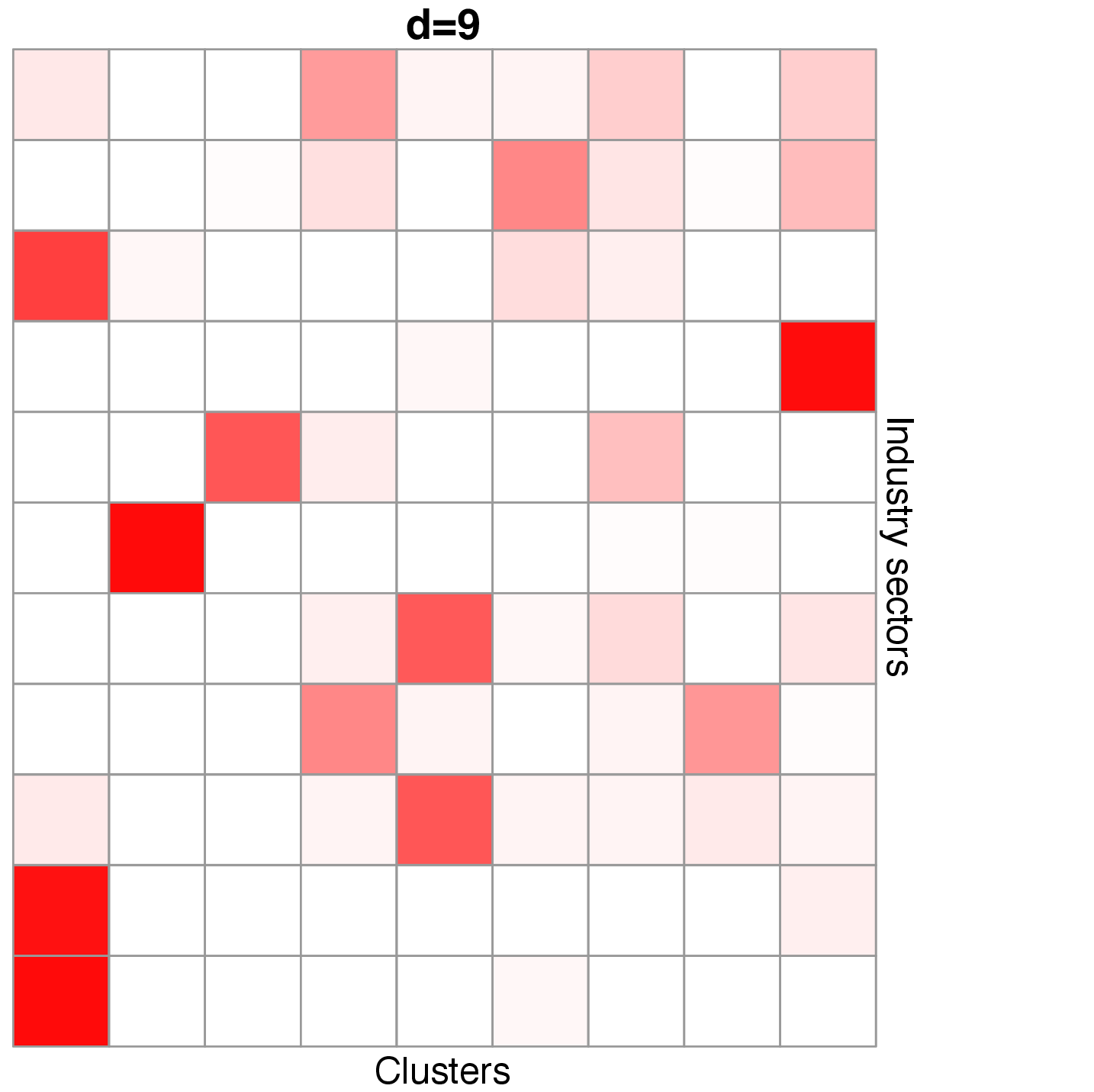, width=5cm, height=6cm} \;\;\;\;
\epsfig{figure=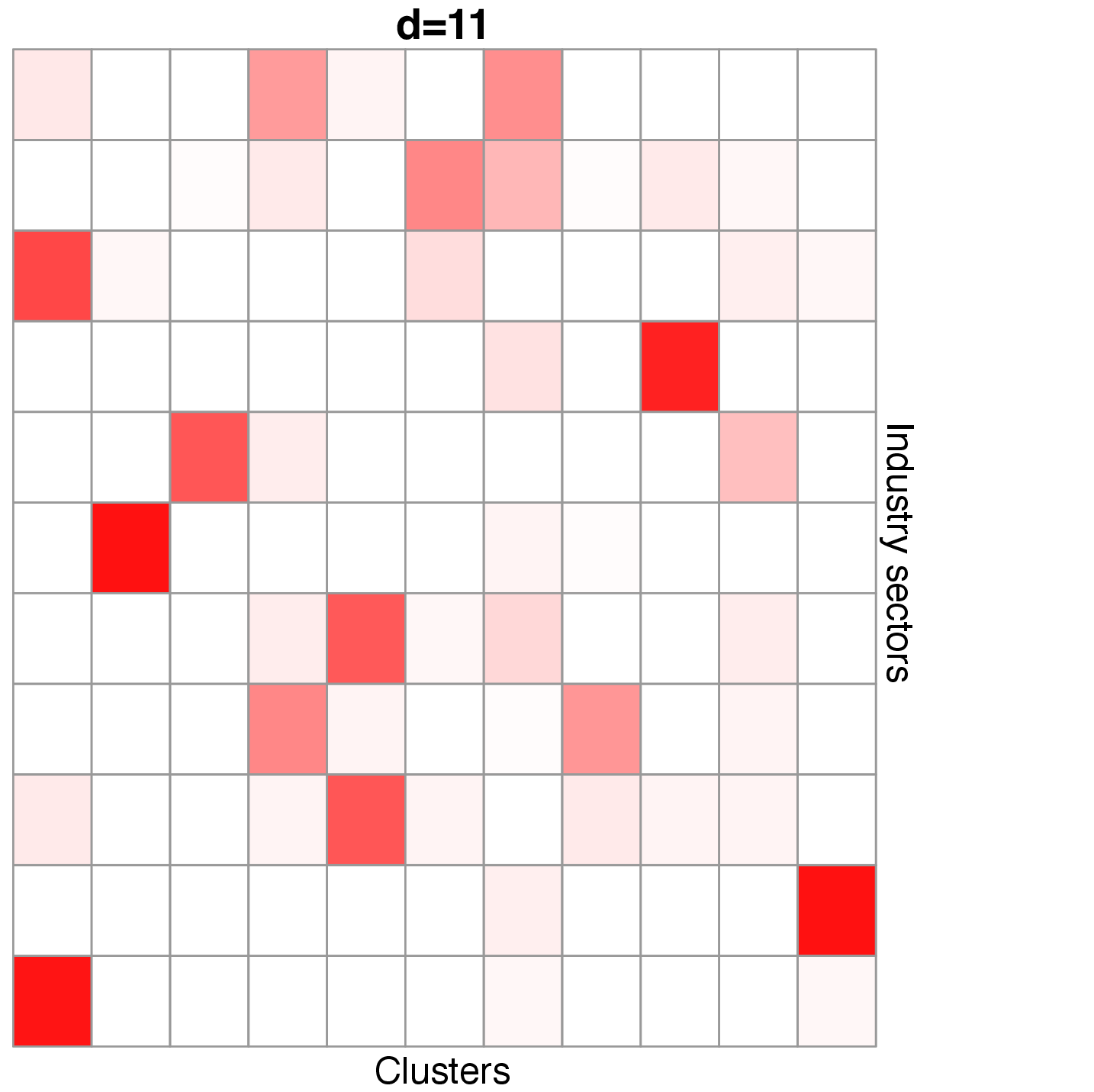, width=6cm, height=6cm}\;\;\;\;
\epsfig{figure=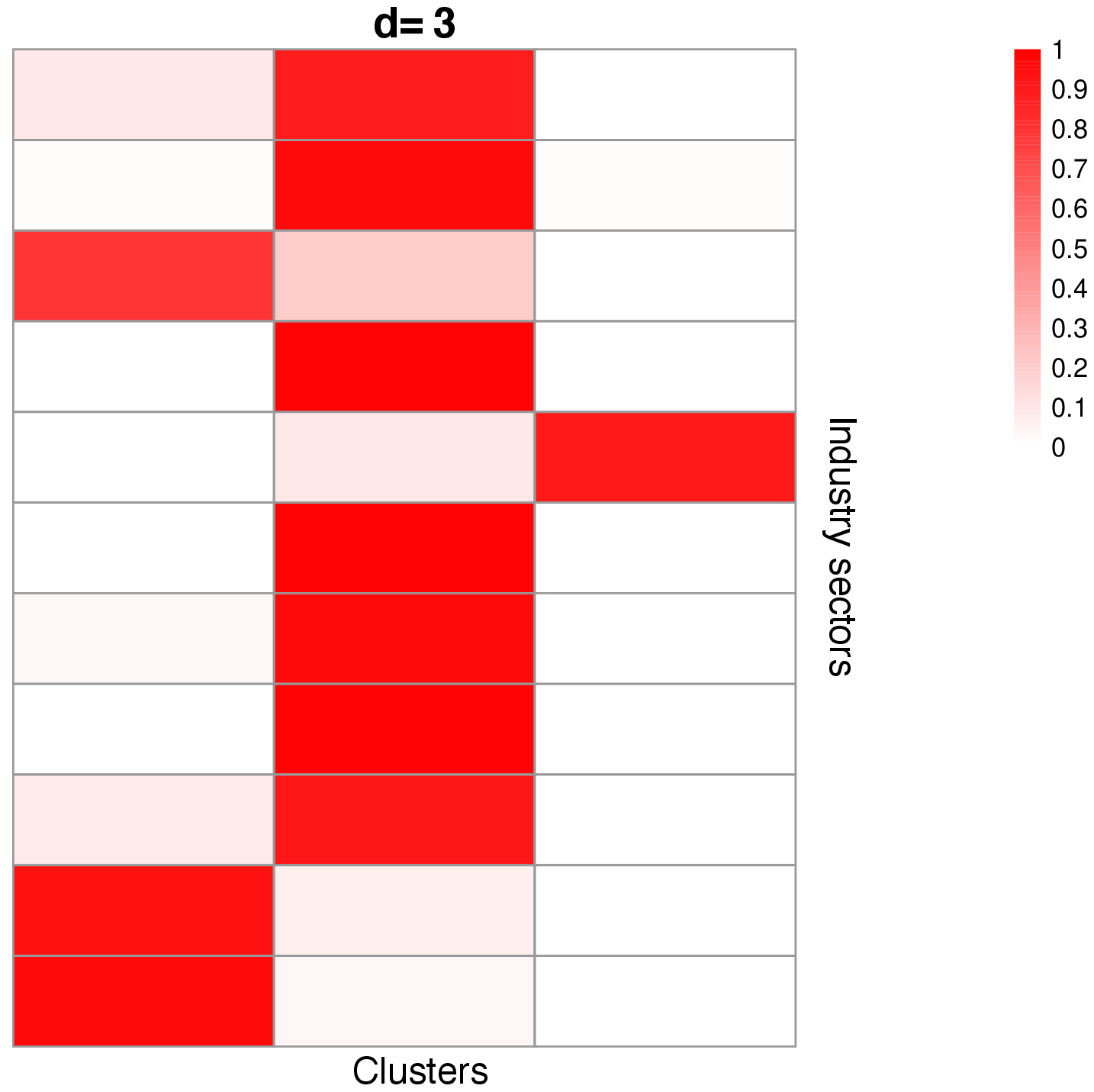, width=4cm, height=6cm}
\caption{\footnotesize{Heat-maps of the distributions of the stocks in each of
the 11 industry sectors (corresponding to 11 rows) over $d$ clusters (corresponding
to $d$ columns), with $d=9, 11 $ and 3. The estimated
numbers of the common and cluster-specific factors are,
respectively, $\wh r_0=1$ and $\wh r =9$.}} \label{fig5}
\epsfig{figure=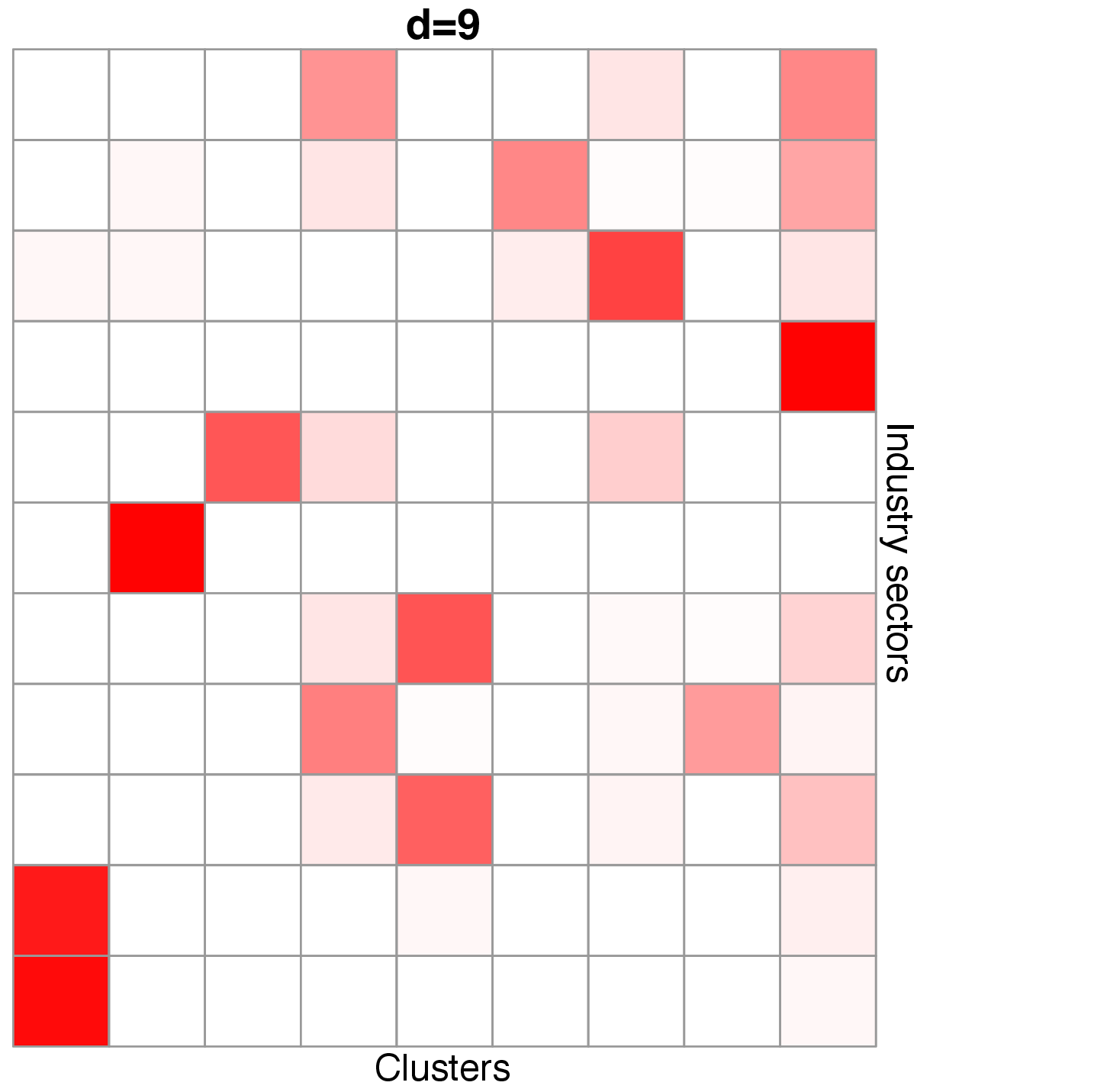, width=5cm, height=6cm} \;\;\;\;
\epsfig{figure=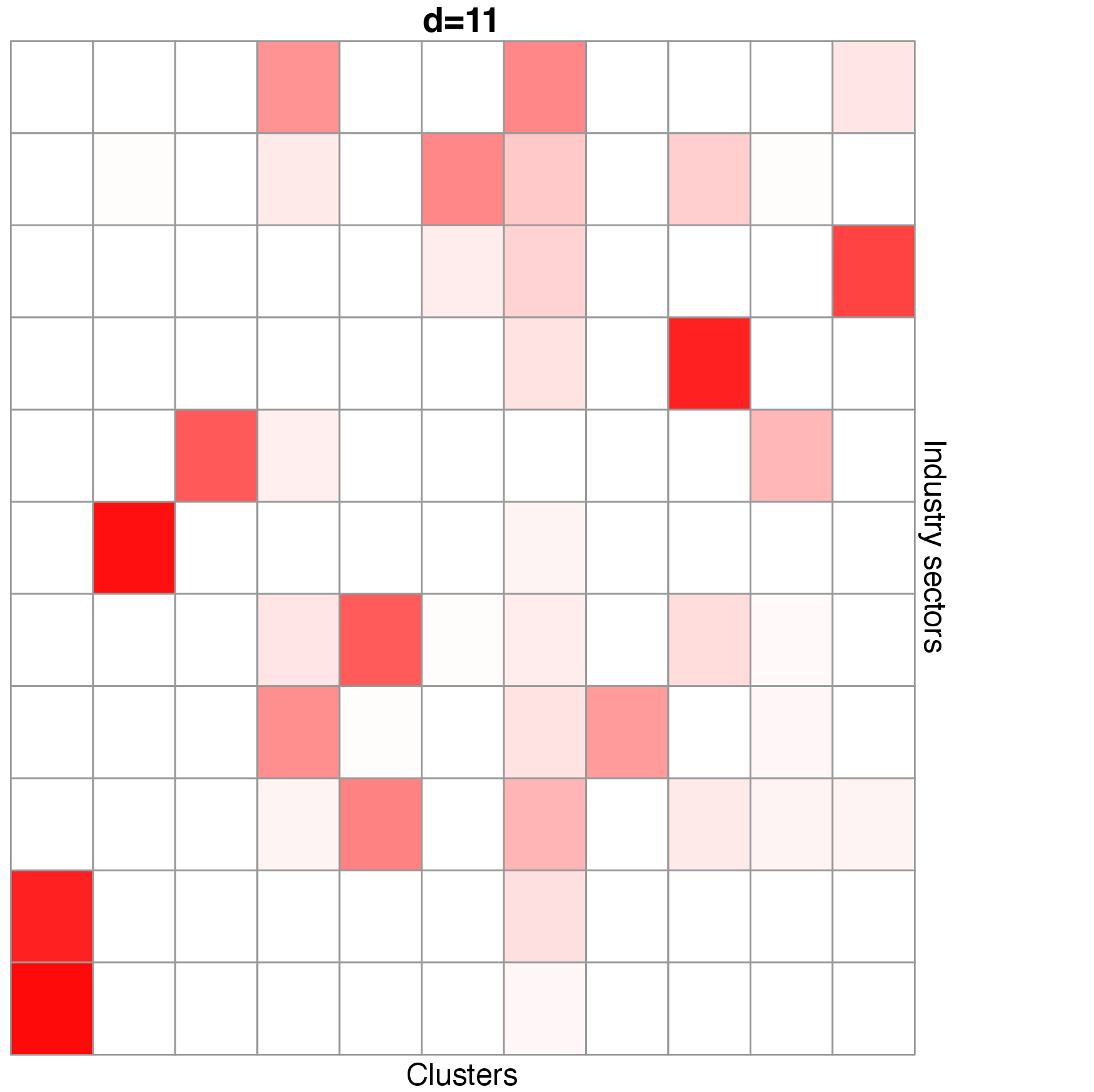, width=6cm, height=6cm}\;\;\;\;
\epsfig{figure=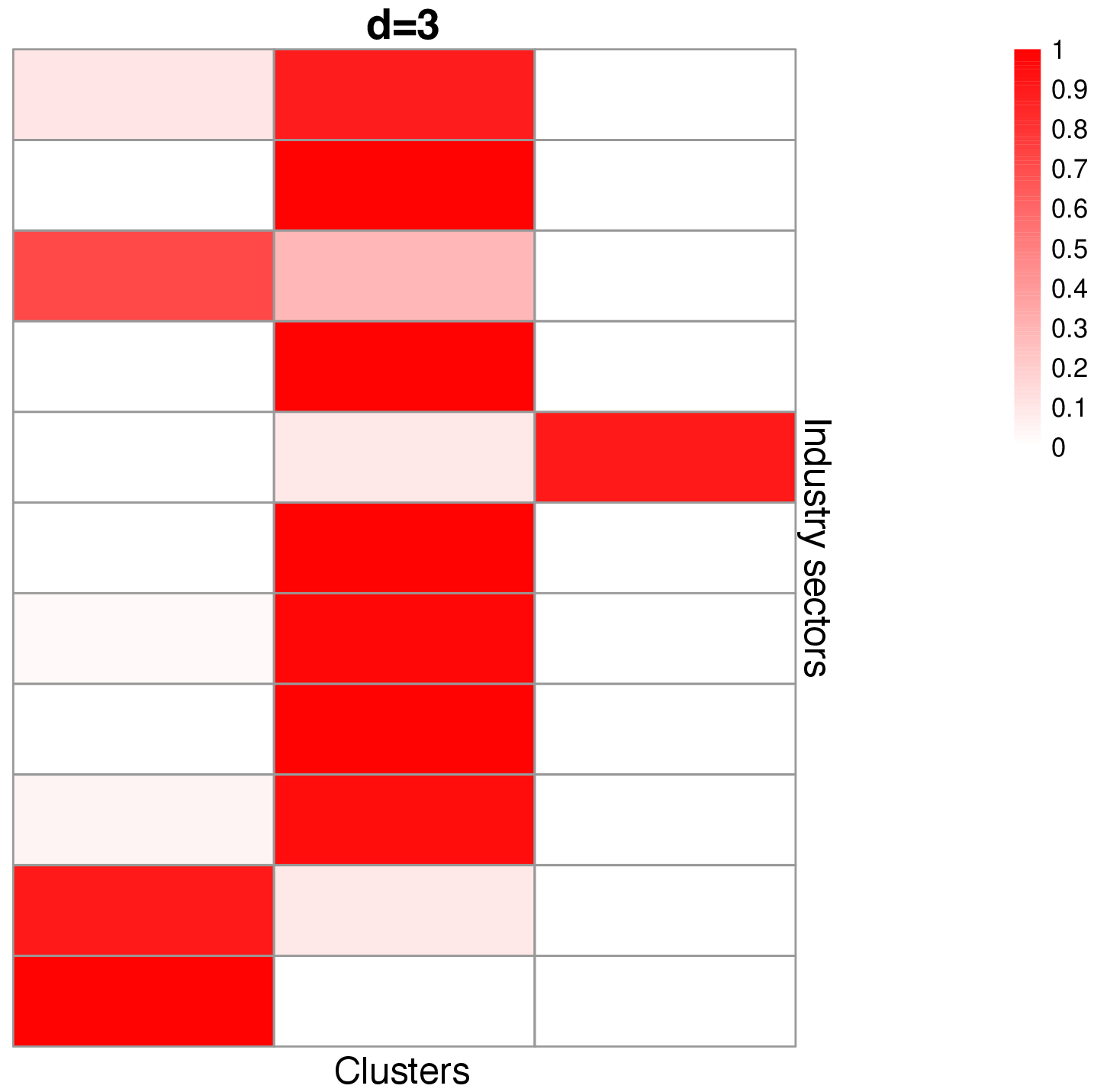, width=4cm, height=6cm}
\caption{\footnotesize{Heat-maps of the distributions of the stocks in each of
the 11 industry sectors (corresponding to 11 rows) over $d$ clusters (corresponding
to $d$ columns), with $d=9, 11 $ and 3. The estimated
numbers of the common and cluster-specific factors are,
respectively, $\wh r_0=1$ and $\wh r =12$.}} \label{fig6}
\end{figure}

If we take  $\wh r_0=1$ and $\wh r=9$, $|\mathfrak{\hat{J}}_{d+1}|=12$ and $\wh p_0= 477-12 = 465$.
The clustering results with $d= 9, 11 $ and 3 are
presented in is presented in Figure \ref{fig5}.
The first cluster mainly consists of the companies in Consumer Staples,
Real Estate and Utilities, Clusters 2 and 3 contain the companies in, respectively,
Health Care and Financials, Cluster 4 contains mainly some companies in Communication
Service and Information Technology, Cluster 5 consists of the companies in Industrials
and Materials, Cluster 6 are mainly the companies in Consumer Discretionary,
Cluster 7 is a mixture of a small number of companies from  each of 5 or
6 different sectors, Cluster 8 is mainly the companies from Information
Technology, Cluster 9 contains almost all companies in Energy.
To examine how stable the clustering is, we also include the results for
$d=11$ and $d=3$ in Figure \ref{fig5}.
When $d$ is increased from 9 to 11, the original Cluster 1 is divided into
new Clusters 1 and 11 with the former consisting of Consumer Staples and
Utilities sectors, and the latter being Real Estate sector. Furthermore the
original Cluster 7 splits into new Clusters 7 and 10, while the other
7 original clusters are hardly changed.
With $d=3$, most companies in each of the 11 sectors stay in one cluster.

\begin{figure}
\centering
\epsfig{figure=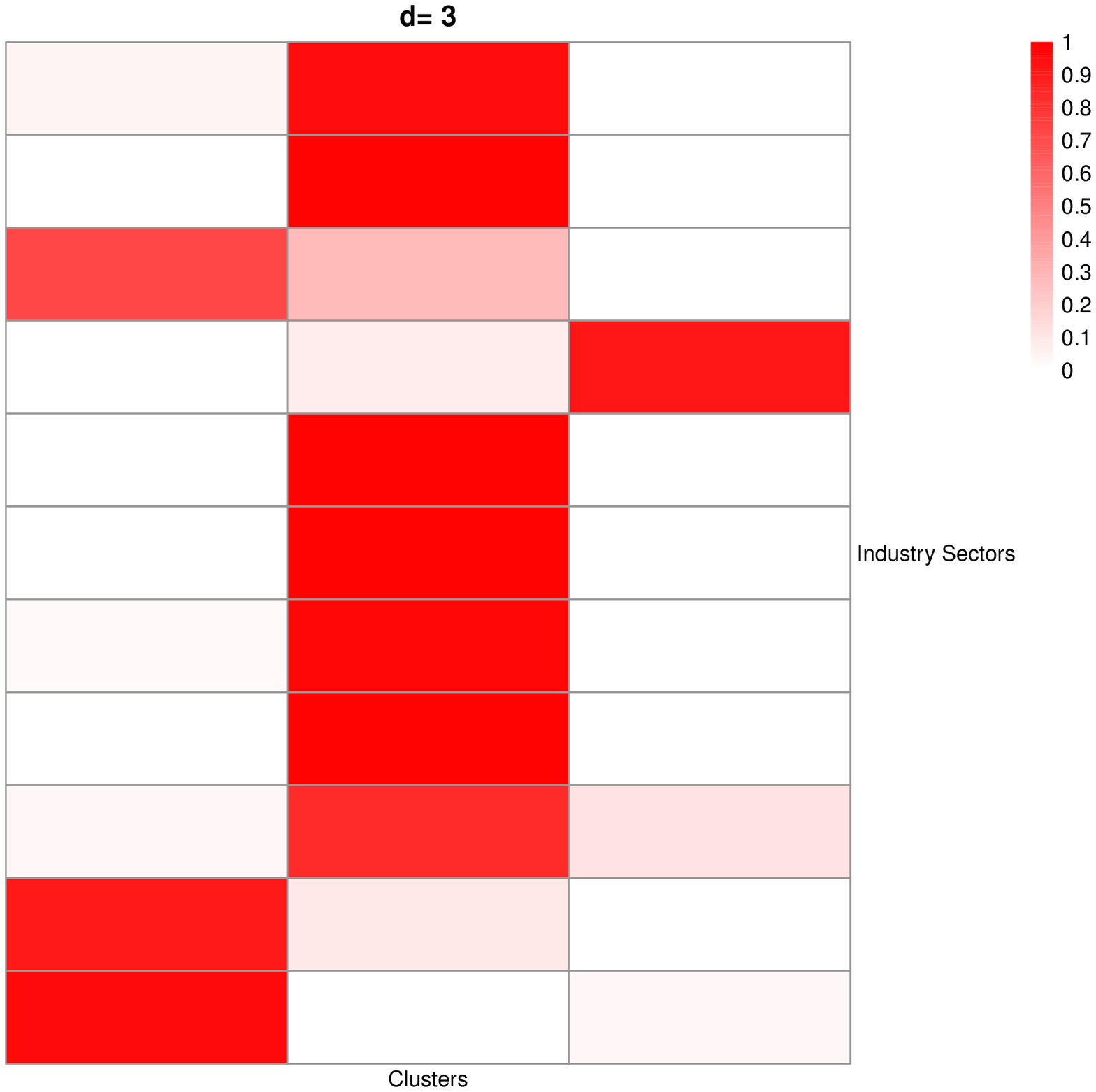, width=5cm, height=6cm}
\epsfig{figure=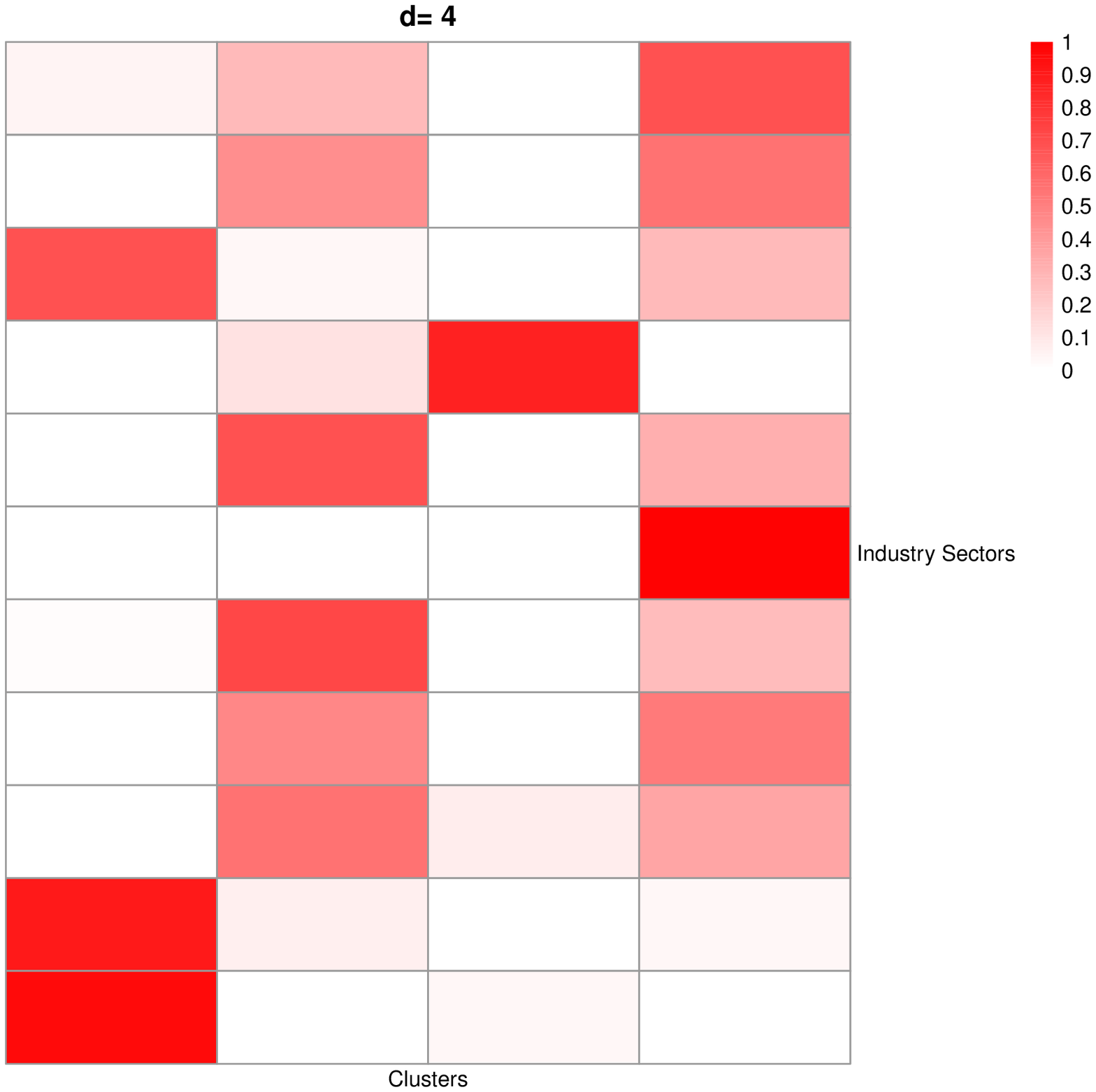, width=5cm, height=6cm}
\caption{\footnotesize{Heat-maps of the distributions of the stocks in each of
the 11 industry sectors (corresponding to 11 rows) over $d$ clusters (corresponding
to $d$ columns) based on Ando and Bai (2017).}} \label{fig5andobai}
\end{figure}

If we take $\wh r_0=1$ and $\wh r=12$, the estimated $\mathfrak{\hat{J}}_{d+1}$
is unchanged. 
The clustering results with $d= 9, 11 $ and 3 are
presented in Figure \ref{fig6}. Comparing with Figure \ref{fig5}, there are some
striking similarities: First the clustering with $d=3$ are almost identical.
For $d=9$, the profiles of Clusters 2, $\cdots$, 6, 8 and 9 are not significantly
changed while Clusters 1 and 7 in Figure \ref{fig5} are somehow mixed together
in Figure \ref{fig6}. With $d=11$, the profiles of Clusters 2 -- 6, 8 -- 10 in
the two figures are about the same while Clusters 7 and 11 are mixed up across
the two figures.

The analysis above indicates that the companies in the same industry sector tend
to share similar dynamic structure in the sense that they are driven by the
same cluster-specific factors. Our analysis is reasonably stable as most the clusters
do not change substantially when the number of the weaker factors chooses different values
 $\wh r= 9$, $\wh r= 12$ or $\wh r= 15$.

 We also apply method of Ando and Bai (2017) to this data set;
 leading to the same estimate $\hat{r}_0=1$, but smaller estimates $\hat{r}=4$ and $\hat{d}=4$.
  The clustering results for  $\hat{d}=3$ and  $\hat{d}=4$ are presented in Figure \ref{fig5andobai} which
  similar to the right parts in Figures \ref{fig5r15}-\ref{fig6}, though the method of Ando and Bai (2017) puts energy companies
  as a separate group. In contrast, our method puts financial companies as a separate group. Note that classical papers in Finance (e.g. Berger and Ofek (1995), Denis et al. (2002), Lemmon et al. (2008)) often eliminate  financial companies from other companies.

\section{Miscellaneous comments}\label{miscomment}
{\bf Robustness}. We identify and distinguish common factors
and cluster-specific factors by the different factor strengths,
i.e. common factors are strong with $\delta =0$, and cluster-specific
factors are weak with $\delta >0$.
However if, for example, one of the common factors has the same
strength as the cluster-specific factors, the number of the strong
factor is then $r_0 -1$ and the number of the weak factors is $r+1$. In this case, the estimated $r_0,\, r$ and the factor loading spaces will all be wrong. Nevertheless a common factor has non-zero loadings on the most components of $\by_t$, hence those loadings must be extremely small in order to be a weak factor.
Therefore its impact on the estimation of the number of clusters, and
the misclassification rates is minor.
%
Simulation results in Supplementary Material  support this assertion.

 \noindent
 {\bf Heterogenous factor strength}.
 We assume two factor strengths: $p$ for common
factors and $p^{1-\delta}$ for cluster-specific factors. If there are $s$
different strengths $p^{1-\delta_1},\cdots,p^{1-\delta_s}$ among
the cluster-specific factors, we can
search for the $s$ largest local maximums among $\wh R_1, \cdots, \wh
R_{J_0-1}$ in Step 1. Moreover, Step 2 should be repeated $s$ times
to estimate $s$ factor loadings corresponding to
different factor strengths. While the asymptotic results can be extended
accordingly,  small values such as $s\le 3$ are sufficient for most practical
applications.

  \noindent{\bf Estimation for $r_j$ and $\bB_j$}.
  When we cluster the time series correctly in Steps 4-5, we obtain the estimator for $\bB_j$ from $\mathbf{\hat{B}}$ directly. We can also run Step 1 on each cluster to estimate  $r_j$. Theorem \ref{thma1k}  ensures the consistency of those
  estimates. 
  Although Theorem  \ref{consistkmean} ensures that most of time series can be clustered correctly, there may be a cluster obtained in Step 5 which is not accuracy enough. It remains an
  open problem to evaluate how the clustering error is propagated
  into the estimation for $r_j$ and $\bB_j$.

\section*{Supplementary Material}
All the technical proofs are presented in an online  supplementary which also contains additional simulation results.
\section*{References}
\begin{description}
\vspace{-3mm}
\item
Aghabozorgi, S., Shirkhorshid, A.S. and Wah, T.Y. (2015).
Time-series clustering -- A decade review. {\sl Information System}, {\bf 53},
16-38.
\vspace{-3mm}
\item
Alonso, A.S. and Pe\~na, D. (2019). Clustering time series by linear dependency.
{\sl Statistics and Computing}. {\bf 29}, 655-676.
\vspace{-3mm}
\item
Ando, T. and Bai, J. (2017). Clustering huge number of financial time series:
a panel data approach with high-dimensional predictors and
factor structures. \JASA, {\bf 519}, 1182-1198.
\vspace{-3mm}
\item
Berger, P.G. and Ofek, E.(1995). Diversification's effect on firm value. Journal of Financial Economics, {\bf 37(1)}, 39-65.

\vspace{-3mm}
\item
Chamberlain, G. (1983). Funds, factors, and diversification in arbitrage pricing models.
{\sl Econometrica}, {\bf 51}, 1305-1323.
\vspace{-3mm}
\item
Chamberlain, G. and Rothschild, M. (1983). Arbitrage, factor structure, and mean-variance
analysis on large asset markets. {\sl Econometrica}, {\bf 51}, 1281-1304.
\vspace{-3mm}
\item
Chang, J., Gao, B. and Yao, Q. (2015). High dimensional stochastic
regression with latent factors, endogeneity and nonlinearity. \JE, {\bf 189},
297-312.
\vspace{-3mm}
\item
 Denis, D. J., Denis, D. K. and Yost, K. (2002). Global Diversification, Industrial Diversification, and Firm Value,
  Journal of Finance, {\bf 57(5)}, 1951-1979.

\vspace{-3mm}
\item
Esling, P. and Agon, C. (2012). Time-series data mining. {\sl ACM Computing
Survey}, {\bf 45}. Article 12.

\vspace{-3mm}
\item
Forni, M., Hallin, M., Lippi, M. and Reichlin, L. (2005).
The generalized dynamic-factor model: one-sided estimation and
forecasting. \JASA, {\bf 100}, 830-840.
\vspace{-3mm}
\item
Fr\"uhwirth-Schnatter, S. and Kaufmann, S. (2008). Model-based clustering of
multiple time series. \JBES, {\bf 26}, 78-89.
\vspace{-3mm}
\item
Hallin, M. and Lippi, M. (2013). Factor models in high-dimensional time series -- a
time-domain approach. \SPA, {\bf 123}, 2678-2695.
\vspace{-3mm}
\item
Kakizawa, Y., Shumway, R.H. and Taniguchi, M. (1998).
Discrimination and clustering for multivariate time series.
\JASA, {\bf 93}, 328-340.
\vspace{-3mm}
\item
Keogh, E. and Lin, J. (2005). Clustering of time-series subsequences is meaningless:
implications for previous and future research.
{\sl Knowledge and Information Systems}, {\bf 8}, 154-177.
\vspace{-3mm}
\item
 Keogh, E. and Ratanamahatana, C.A. (2005). Exact indexing of dynamic time warping.
{\sl Knowledge and Information Systems}, {\bf 7}, 358-386.
\vspace{-3mm}
\item
Khaleghi, A., Ryabko, D., Mary, J. and Preux, P. (2016).
Consistent algorithms for clustering time series.
\JMLR, {\bf 17}, 1-32.
\vspace{-3mm}
\item
Lam, C. and Yao, Q. (2012). Factor modelling for high-dimensional time series: inference
for the number of factors.  \AS, {\bf 40}, 694-726.

\vspace{-3mm}
\item
Lemmon, M. L., Roberts, M.R. and Zender, J. F. (2008). Back to the Beginning: Persistence and the
Cross-Section of Corporate Capital Structure. Journal of Finance, {\bf 63(4)}, 1575-1608.

\vspace{-3mm}
\item
Li, Z., Wang, Q. and Yao, J. (2017). Identifying the number of factors from
singular values of a large sample auto-covariance matrix. \AS, {\bf 45},
257-288.
\vspace{-3mm}
\item
Liao, T.W. (2005). Clustering of time series data -- a survey.
{\sl Pattern Recognition}, {\bf 38}, 1857-1874.

\vspace{-3mm}
\item
Maharaj, E.A., D'Urso, P. and
Caiado, J. (2019). {\sl Time Series Clustering and Classification}.
Chapman and Hall/CRC.

\vspace{-3mm}
\item
Pe{\~n}a, D. and Box, E.P. (1987). Identifying a simplifying structure in time
series. \JASA, {\bf 82}, 836-843.
\vspace{-3mm}
\item
Pe{\~n}a, D. and Poncela, P. (2006). Nonstationary dynamic factor analysis. \JSPI, {\bf 136}, 1237-1257.
\vspace{-3mm}
\item
Roelofsen, P. (2018). Time series clustering. Vrije Universiteit Ansterdam.\\
{\tt https://www.math.vu.nl/$\sim$sbhulai/papers/thesis-roelofsen.pdf}.

\vspace{-3mm}
\item
Vershynin, R. (2010). Introduction to the non-asymptotic analysis of random matrices. \\
{\tt arXiv.1011.3027}.
\vspace{-3mm}
\item
Yao, Q., Tong, H., Finkenst\"adt, B. and Stenseth, N.C. (2000).
Common structure in panels of short ecological time series.
        {\sl  Proceeding of the Royal Society (London)}, {\bf B}, {\bf 267}, 2457-2467.

\vspace{-3mm}
\item
Zhang, T. (2013). Clustering high-dimensional time series based on parallelism.
\JASA, {\bf 108}, 577-588.
\vspace{-3mm}
\item
Zolhavarieh, S., Aghabozorgi, S. and Teh, Y.W. (2014). A
review of subsequence time series clustering. {\sl The Scientific World
Journal}, Article 312512.

\end{description}
%
%
%
\end{document}